\documentclass[12pt]{amsart}
\usepackage[left=1.5cm]{geometry}
\usepackage[all]{xy}
\usepackage{stmaryrd, amssymb}
\usepackage{amscd,verbatim}
  
\newtheorem{theorem}{Theorem}[section]
\newtheorem{proposition}[theorem]{Proposition}
\newtheorem{corollary}[theorem]{Corollary}
\newtheorem{definition}[theorem]{Definition}

\newtheorem{lemma}[theorem]{Lemma}
\newtheorem{remark}[theorem]{Remark}

\newcommand{\Z}{\mathbb{Z}}
\newcommand{\C}{\mathcal{C}}
\newcommand{\D}{\mathcal{D}}

\newcommand{\B}{\mathcal{B}}
\newcommand{\F}{\mathbb{F}}
\newcommand{\G}{\mathbb{G}}

\newcommand{\q}{\mathcal{P}}
\newcommand{\Q}{\mathbb{Q}}

\newcommand{\Spl}{\mathbf{Spl}}

\newcommand{\Set}{\operatorname{\mathbf{Set}}}
\newcommand{\Ord}{\operatorname{\mathbf{Ord}}}
\newcommand{\Top}{\operatorname{\mathbf{Top}}}
\newcommand{\Simpl}{\mathbf{Simpl}}
\newcommand{\Dim}{{\operatorname{dim}}}
\newcommand{\tf}{{\operatorname{tf}}}
\newcommand{\coh}{{\operatorname{coh}}}
\newcommand{\Ind}{{\operatorname{Ind}}}

\renewcommand{\epsilon}{\varepsilon}

\title{Algebraic K-Theory and Modular Symbols}
\author{Fei SUN}

\begin{document}
\sf

\thanks{fei.sun@icloud.com}
\maketitle

\include{resume}

\section*{Abstract}
In this paper, we calculate the differential $d^1$ of the rank spectral sequence defined in \cite{K}. We generalize Quillen's spectral sequence given in \cite{Q2} from Dedekind domain to general integral Noetherian ring $A$ by considering the Q-construction $Q^{\tf}(A)$ of the category of finitely generated torsion-free modules. In particular, by resolution theorem, if $A$ is regular, then $Q^{\tf}(A)$ is homotopy equivalent to $QP(A)$ where $P(A)$ is the category of finitely generated projective $A$-modules.  We deduce the differential $d^1$ by using Tits buildings, Steinberg modules, and modular symbols in the sense of Ash-Rudolph defined in \cite{AR}. The spirit of Quillen's categorical homotopy theory will be used intensively throughout this paper.

\setcounter{tocdepth}{4} 

\tableofcontents  


\section*{Introduction}
\addcontentsline{toc}{chapter}{Introduction}
\subsection*{Bird's-Eye View of the Paper}
\addcontentsline{toc}{section}{Bird's-Eye View of the Paper}
In his preprint \cite{K}, Bruno Kahn introduced a notion of cellular functor (\cite[definition 2.3.2]{K}), that is, a functor $T:\C\to\D$ in the category of small categories $\mathbf{Cat}$ which is fully faithful and such that $Hom_{\D}(d, T(c))=\emptyset$ for any $c\in\C$ and $d\in\D-T(\C)$. We have the covariant functor 
$$\F_T: \D\to\mathbf{Cat};\ \ \ \ \ d\mapsto T\downarrow d.$$
Here, $T\downarrow d$ is the comma category (see Section 1.1 for its definition).
We denote $N:\mathbf{Cat}\to \mathbf{sSet}$ the nerve functor such that for any small category $\C$ we have $N\C\in\mathbf{sSet}$ with
$$
N_i\C=\coprod_{c_0\to\cdots\to c_i}*,\ \ \ \ \ \ c_i\in\C,\ \ i\geq 0.
$$
Here, $*$ is the singleton.  
If $X$ is a simplicial abelian group, we denote by $C_*(X):=[X]$ the singular chain complex associated to $X$ such that differentials are taken to be the alternating sum of faces. If $X$ is a simplicial set, we will write $\mathbb{Z}X$ for the free simplicial abelian group such that $(\mathbb{Z}X)_n=\mathbb{Z}X_n$ and $C_*(X):=[\Z X]$. For a small category $\C$, we write $C_*(\C):=C_*(N\C)=[\Z N\C]$.
Moreover, if $X$ is a bisimplicial set, then we may take the associated diagonal simplicial set $\delta X$ with
$$
\delta_nX=X_{n,n}
$$ 
and $C_*(X):=[\Z\delta X]$. 

\begin{definition}\label{singchain}
If $\F: \D\to\mathbf{sSet}$ is a functor taking value in the category of simplicial sets, then we define $N(\D,\F)$ as a bisimplicial set with
$$
N_{p,q}(\D,\F):=\coprod_{d_0\to\cdots\to d_p}\F_q(d_0),
$$
and $C_*(\D,\F):=[\Z\delta N(\D,\F)]$.
\end{definition}
In particular, for $\F_T:\D\to\mathbf{Cat}$ given above and $N\F_T: \D\to\mathbf{sSet}$,
we have:

\begin{enumerate}
\item $N(\D,\F_T)$ to be the bisimplicial set with
          $$
          N_{p,q}(\D,\F_T):=\coprod_{d_0\to\cdots\to d_p} N_q(\F_T(d_0))=\coprod_{d_0\to\cdots\to d_p}\coprod_{c_0\to\cdots\to c_q\to d_0}*
          $$
      such that $d_i\in\D$ for $0\leq i\leq p$ and $c_j\in\C$ for $0\leq j\leq q$.
          
\item
$C_*(\D,\F_T)=[\Z \delta N(\D,\F_T)]$ to be the singular chain complex of the diagonal simplicial set $\delta N(\D,\F_T)$.
\end{enumerate}

We notice that this construction in definition \ref{singchain} makes $C_*(\D, \F)$ functorial with respect to $\F$. Let $\star$ be the constant functor taking value the terminal object (singleton) $*\in\mathbf{Cat}$ with one object and one morphism, we define $C_*(\D,\widetilde{\F_T})$ to be the homotopy fiber of
$C_*(\D,\F_T)\to C_*(\D,\star)$
in the derived category of Abelian groups $\mathbf{D}(Ab)$ so that it gives rise to a distinguished triangle
$$
C_*(\D,\widetilde{\F_T})\to C_*(\D,\F_T)\to C_*(\D)\to C_*(\D,\widetilde{\F_T})[1].
$$
Kahn constructed in \cite[proposition 2.3.5]{K} a naturally commutative diagram of categories
$$
\xymatrix{
(\D-\C)\int\F_T\ar[r]\ar[d] & \C\ar[d]\\
\D-\C\ar[r] & \D
}
$$
and proved that it becomes homotopy cocartesian after applying the nerve functor. Combined with Thomason's theorem (theorem \ref{thomason}), we get a long exact sequence (c.f, \cite[Thm 2.3.6]{K})
$$
\cdots\to H_i(\D-\C,\widetilde{\F_T})\to H_i(\C)\to H_i(\D)\to H_{i-1}(\D-\C, \widetilde{\F_T})\to\cdots
$$
\\

This technique can be applied to the case of Quillen's higher K-theory. 
Suppose that $A$ is an integral Noetherian domain and $Q^{\tf}(A)$ be Quillen's Q-construction over the the category of finitely generated torsion-free modules. For an $A$-module $M$ we define its rank as 
$$
rank(A):=\dim(M\otimes_AK)
$$ 
where $K=Frac(A)$. Then we denote $Q^{\tf}_n(A)$ for $n\geq 0$ the full subcategory of $Q^{\tf}(A)$ of modules with rank less or equal to $n$. 
If $d$ is a torsion-free Noetherian $A$-module, a submodule $c\subset d$ is said to be pure if $d/c$ is torsion-free. According to \cite[Prop. 4.2.4]{K}, there is a bijection between the poset of proper pure submodules of $d$ and the poset of proper subspaces of $V:=d\otimes K$. 
Let  $V$ be an $n$-dimensional vector space over a field $K$. According to \cite{Q2}, the Tits building $T(V)$ is defined to be the simplicial complex whose $p$-simplices are
$$V_0<\cdots <V_p;\ \ \ \ \ V_0>0,\ V_p<V.$$
We use Quillen's model of the suspension of the Tits building, $\Sigma T(V)$, which is a simplicial complex where a $p$-simplex is a flag of subspaces
$$
W_0<\cdots<W_p
$$
with either $0<W_0$ or $W_p<V$. This does give a model of suspension of $T(V)$ for $\dim(V)>1$, and when $\dim(V)=1$, $\Sigma T(V)$ consists of two points, whereas $T(V)=\emptyset$.
We denote by $J(V)$ the ordered set of proper layers in $V$. This consists of pairs $(W_0,W_1)$ of subspaces of $V$, with $W_0\subset W_1$, excluding the pair $(0,V)$, and where $(W_0,W_1)\leq (W'_0,W'_1)$ if $W'_0\subset W_0\subset W_1\subset W'_1$. $J(V)$ has the homotopy type of the suspension of the Tits building, $\Sigma T(V)$.
Then by \cite[Cor. 4.2.6]{K}, the functor $Q^{\tf}_{n-1}(A)\downarrow d\to J(V)$ is an equivalence of categories.  
\\

The Solomon-Tits theorem (\cite[Theorem 2]{Q2}) says that $T(V)$ has the homotopy type of a bouquet of $(n-2)$-spheres for $n\geq 2$.
The Steinberg module $St(V)$ will be defined as 
$$
St(V):=\left\{
\begin{array}{ccccc}
H_{n-2}(T(V)) & ; & n\geq 2\\
\Z & ; & n=1
\end{array}
\right.
$$
\begin{definition}\label{redms}[\cite{KS}, (1)]
The reduced Steinberg module, $\widetilde{St}(V)$, is defined as
$$
\widetilde{St}(V):=\left\{
\begin{array}{ccc}
St(V) & ;&  \Dim(V)>2\\
Ker(St(V)\to\Z) & ; & \Dim(V)=2\\
\Z & ; & \Dim(V)=1\\
\Z & ; & \Dim(V)=0
\end{array}
\right.
$$
\end{definition}
\

Let us write $Q_{n-1}:=Q^{\tf}_{n-1}$, $Q_n:=Q^{\tf}_n$ for short and denote $T_n: Q_{n-1}\hookrightarrow Q_n$ the inclusion functor (which is connected cellular by definition).
By \cite[Thm 2.3.6]{K}, we are able to describe the homotopy mapping cone of $C_*(Q_{n-1})\hookrightarrow C_*(Q_n)$ as $C_*(Q_n-Q_{n-1}, \widetilde{\F_{T_n}})[1]$.
If we denote $D^1_{p,q}=H_{p+q}(Q_p)$ and $E^1_{p,q}=H_{p+q-1}(Q_p, \widetilde{\F_{T_p}})$, we
get an exact couple and hence a spectral sequence (\cite[theorem 2.4.1]{K}) converging to the homology groups of $BQ^{\tf}(A)$
$$
E^1_{p,q}=H_{p+q-1}(Q_p(A)-Q_{p-1}(A), \widetilde{\F_p})\Rightarrow H_{p+q}(BQ^{\tf}(A)).
$$
This spectral sequence is called the rank spectral sequence. 
Moreover, since for any $d\in Q_n-Q_{n-1}$, $\F_{T_n}(d)=T_n\downarrow d\simeq J(V)$ (c.f, Lemma \ref{layer}), the rank spectral sequence becomes
$$
E^1_{n,i-n+1}=H_i(Q_n-Q_{n-1}, \widetilde{\F_{T_n}})\simeq 
\bigoplus_{d}H_{i-n+1}(Aut(d), \widetilde{St}(V))
\Rightarrow H_{i+1}(BQ^{\tf}(A))
$$
where $d$ runs over the isomorphism classes of torsion-free $A$-modules of rank $n$.
We denote by $QM(A)$ (resp. $QP(A)$) the Quillen's Q-construction over the category of finitely generated $A$-modules (resp. finitely generated projective $A$-modules).
By Quillen's resolution theorem, the inclusion $Q^{\tf}(A)\hookrightarrow QM(A)$ is a homotopy equivalence (c.f, \cite[4.2.3]{K} ). On the other hand, if $A$ is regular, by \cite[Corollary 2, p26]{Q1}, the inclusion $QP(A)\hookrightarrow QM(A)$ is a homotopy equivalence. So in the case of $A$ being regular, the rank spectral sequence converges to $H_*(BQP(A))$.
\\

We notice that, in particular, a Dedekind domain is regular. So our construction above generalizes Quillen's setup in \cite[Section 3]{Q2}.
The idea above applies to the case of an integral scheme $X$ and the Q-construction over the category of finitely generated coherent sheaves $Q^{\coh}(X)$ (resp. finitely generated torsion-free sheaves $Q^{\tf}(X)$). Please go to \cite[Section 4]{K} for more details.
\\

To write down the formula of $d^1$, we use Ash-Rudolph's definition of modular symbols $[g_1,\cdots,g_n]$ with $g_1,\cdots,g_n$ non-zero vectors of $V$. 
However, their definition takes the homology of Tits buildings rather than reduced homologies. 
So there is a slight mistake in the paper \cite{AR} in the case of $n=\Dim(V)=2$. 
This mistake is corrected in a joint paper of Bruno Kahn and Fei Sun in \cite{KS} by taking reduced Steinberg modules (Definition \ref{redms}).
We also notice that our construction
of rank spectral sequence always use \textit{suspension} and \textit{reduced} homologies.
Under this perspective,  we define a new symbol, called extended symbols, $\partial(0,g_1,\cdots,g_n)$ with $g_1,\cdots,g_n\in V-\{0\}$ that form a set of generators of $\widetilde{St}(V)$. In section 4.4, we show that modulo certain isomorphism we may identify the symbols $\partial(0,g_1,\cdots,g_n)$ with $[g_1,\cdots,g_n]$ for $n\geq 2$. 
\\

In Chapter 3, we deduce that $d^1$ of rank spectral sequence can be calculated via distinguished triangles of coefficients/functors. 
By \cite[Prop. 4.2.4]{K}, we may use alternatively the notations $\widetilde{St}(d)=\widetilde{H}_{n-1}(Q_n^{\tf}(A)\downarrow d)$ and $\widetilde{St}(V)=\widetilde{H}_{n-1}(J(V))$ since $Q_n^{\tf}(A)\downarrow d$ and $J(V)$ are canonically isomorphic.  
Our calculation by using cellular functors gives us a morphism (c.f, (\ref{coeff}))
$$
\widetilde{St}(d)\to\bigoplus_{c\to d\in\F_{T'}(d)}\widetilde{St}(c).
$$
We divide this morphism into two parts:
\begin{itemize}
\item the direct sums indexed by admissible monomorphisms $c\rightarrowtail d$, which is calculated in section 5.2.1.
\item the direct sums indexed by admissible epimorphisms $d\twoheadrightarrow d$, which is calculated in section 5.2.2.
\end{itemize}
Then in section 5.2.3, we combine the above two parts to get the formula for $d^1$ on coefficients. Finally, in Chapter 6, we apply the homology functor $H_q(\D-\C,-)$ to the formula on coefficients to get the formula for $d^1$.
\\

\subsection*{What Happens Next}
\addcontentsline{toc}{section}{What Happens Next}
Consider the following situation (c.f: \cite{QG}): we are given $C$ a non-singular projective curve over a finite field $k$, and $X=Spec(A)$ is the complement
of a closed point in $C$. 
The next step is to calculate the ranks of $K_i(A)$ for $i\geq 0$ by using $d^1$. The ranks have been obtained by G. Harder (\cite[3.2.3]{H}) via completely different method which asserts that $K_i(A)$ is finite for $i>0$. This paper makes preparation for a \textit{functorial} calculation. 
Although the result is not attained by this paper, it will be a very interesting follow-up project.
\\

Actually, Quillen showed that $A$ satisfies
\begin{enumerate}
\item $Pic(A)$ is finite (this is actually a classical result).
\item If $d$ is a finitely generated projective $A$-module and $V=d\otimes_AK$, then $H_i(Aut(d), St(V))$ is a finitely generated Abelian group for every $i$.
\end{enumerate}
In particular, for $\Gamma=Aut(d)$ and $d\in\D-\C$ we may find a normal subgroup $\Gamma'$ of $\Gamma$ of finite index so that there is a Hochschild-Serre spectral sequence
$$H_p(\Gamma/\Gamma', H_q(\Gamma', \widetilde{St}(d)))\Longrightarrow H_{p+q}(\Gamma, \widetilde{St}(d)).$$
Since the quotient group 
$\Gamma/\Gamma'$ is finite, the homology group $H_p(\Gamma/\Gamma', H_q(\Gamma', \widetilde{St}(d)))$ is of torsion for $p\neq 0$. 
Furthermore, Quillen also showed that $\Gamma'$ can be chosen so that $St(d)$ is a $\Z\Gamma'$-projective module of finite type. Hence the homology groups $H_q(\Gamma', \widetilde{St}(d))$ is finitely generated for $q=0$ and trivial otherwise (\cite{QG}, section 5).
So, the rank spectral sequence satisfies:
$$E^1_{p,q}\otimes\mathbb{Q}=0;\ \ \ q\neq 0$$
Therefore, we have 
\begin{equation}\label{e2}
E^2_{p,0}\otimes\mathbb{Q}=E^{\infty}_{p,0}\otimes\mathbb{Q}=H_p(BQP(A), \mathbb{Q}).
\end{equation}
Since $QP(A)$ is a homotopy commutative and associative H-space, the graded homology group $H_*(QP(A); \Q)$ forms a Hopf algebra with commutative comultiplication, in particular, the Lie bracket $[-,-]$ is zero. We consider its primitive part $P_*:=PH_*(QP(A);\Q)$ which is a sub Lie algebra of $H_*:=H_*(QP(A);\Q)$.  
The universal enveloping algebra $UP_*$ is defined as
$$
UP_*=T(P_*)/I
$$
where $T(P_*)=\bigoplus_nP_*^{\otimes n}$ is the universal tensor algebra and $I$ is the ideal of $TP_*$ generated by $x\otimes y-(-1)^{ij}y\otimes x-[x,y]$ with any $x\in P_i$ and $y\in P_j$. 
According to \cite{MM}, $H_*(QP(A);\Q)$ equals to the universal enveloping algebra of $PH_*(QP(A);\Q)$, that is: 
$$
H_*=H_*(QP(A);\Q)=U(PH_*(QP(A);\Q))=UP_*.
$$
By a theorem of Cartan and Serre (cf. \cite[Appendix]{MM}), the graded homotopy group $K_*(A)=\pi_{*+1}(QP(A);\Q)$ is canonically isomorphic to $P_{*+1}$ as Lie algebras. 
\\

Now, we have known $d^1$. So if one can use it to calculate $E^2_{i,0}$ then we obtain the ranks of
$E^2_{i,0}\otimes\Q=H_i(QP(A); \Q)$. Since $A$ is a Dedekind domain, we have $K_0=\Z\oplus Cl(A)$ and hence $rank(K_0\otimes\Q)=rank(P_1)=1$ as $Cl(A)=Pic(A)$ is finite. Finally, since $H_*=UP_*$, we may be able to calculate the ranks of $K_i(A)$ by induction on degrees $i$.
\\

\subsection*{Rank Spectral Sequence and (Extended) Modular Symbols-An Example}
\addcontentsline{toc}{section}{An Example}
Let $A$ be an integral Noetherian domain.
Recall that if we take $Q_n=Q^{\tf}_n(A)$ for any $n\geq 0$ and $T_n: Q_{n-1}\to Q_n$ the natural inclusion functor (which is connected cellular), then the rank spectral sequence takes the form:
$$
E^1_{n,i-n+1}=\bigoplus_d H_{i-n+1}(Aut(d), \widetilde{St}(V))\Longrightarrow H_{i+1}(BQ^{\tf}(A)).
$$
Let us describe $E^1_{n,i-n+1}$ in a low dimensional case, i.e, $n=1$, using extended symbols. We will see (in Chapter 4) that an advantage of the extended symbols is that they generate the reduced Steinberg module $\widetilde{St}(V)$ even when $n=\dim(V)=1$ (while modular symbols are used only for dimension at least two). Thus, we can use extended symbols to calculate $E^1_{1,q}$. In this case, $Aut(d)\simeq A^{\times}$, hence
$$
E^1_{1,q}\simeq \bigoplus_d H_q(A^{\times}, \widetilde{St}(V))
$$ 
where $d$ runs over the isomorphism classes of torsion-free $A$-modules of rank one.
However, since $\partial(0,g)=\partial(0,ag)=a\cdot\partial(0,g)$ for any $a\in A^{\times}$, the action of $A^{\times}$ on $\widetilde{St}(V)=\Z$ is trivial. Thus we get
$$
E^1_{1,q}\simeq\bigoplus_d H_q(A^{\times}).
$$
In particular, we have 
$$
E^1_{1,0}\simeq\bigoplus_d\Z.
$$
\

Moreover, since modular symbols form a set of generators of reduced Steinberg modules satisfying certain relations, we can use modular symbols to do some calculation. Let us consider $A=\F_q[t]$ with $\F_q$ a finite field such that $q=p^n$ for some prime $p$ and integer $n$. We will apply \cite[Theorem 4.1]{LS} and its Corollary. However, as stated in Introduction, we should replace $St(V)$ with $\widetilde{St}(V)$ in \cite[Theorem 3.1]{LS}. Since $A$ is an Euclidean domain, \cite[Theorem 4.1]{LS} implies that $E^1_{n,0}=0$ for $n\geq 2$.
So we obtain
$$
d^1_{n,0}=0,\ \ \ \ \forall n\geq 2.
$$
Finally, we adopt the notations given in the previous section. By our discussion there, we have
$$
E_{p,0}^2\otimes\Q=E_{p,0}^{\infty}\otimes\Q=H_p(BQP(A),\Q)=0,\ \ \ \ \ \ p\geq 2.
$$
Since $P_i$ are direct summands of $H_i$ and $H_i=0$ for $i\geq 2$, it is immediate that $P_i=0$ for $i\geq 2$. Hence, 
$$
rank(K_i(A)\otimes\Q)=rank(P_{i+1})=
\left\{
\begin{array}{cc}
1; & i=0\\
0; & i>0
\end{array}
\right.
$$
which coincide with the result of G. Harder in \cite{H}.


\subsection*{Notations and Conventions}\

We will denote $\widetilde{H_*}(C_*)$ the reduced homology of the chain complex $C_*$. 
We use the notation $\simeq$\ or\ $\xrightarrow{\simeq}$ to denote an isomorphism and $\sim$\ or\ $\xrightarrow{\sim}$ to denote a quasi-isomorphism. 
We still use the notation $\simeq$\ or\ $\xrightarrow{\simeq}$ to denote a homeomorphism. On the other hand,   the notation $\approx$ or $\xrightarrow{\approx}$ denotes a homotopy equivalence .

\newpage

\section{General Setup}

\subsection{Basic Definitions and Properties}
Let $\C,\ \D$ be two categories and $\mathbf{Cat}$ be the category of small categories. Recall that if $T:\C\to\D$ is a functor, then for any $d\in\D$, the comma category $T\downarrow d$ has objects $T(c)\to d\in Mor(\D)$ and a morphism $(T(c_0)\to d_0)\to(T(c_1)\to d_1)\in Mor(T\downarrow d)$ is given by morphisms $(c_0\to c_1)\in Mor(\C)$ and $d_0\to d_1\in Mor(\D)$ which fit into the commutative diagram
$$
\xymatrix{
T(c_0)\ar[r]\ar[d] & d_0\ar[d]\\
T(c_1)\ar[r] & d_1
}
$$ 

\begin{definition}\label{cell}
Let $T: \C\to\D$ be a functor, we call $T$ a cellular functor if it satisfies
\begin{description}
\item[(1)] $T$ is fully faithful;
\item[(2)] $Hom_{\D}(d, T(c))=\emptyset$ for any $c\in\C$ and $d\in \D-T(\C)$
\end{description}
Moreover, we say that a cellular functor $T: \C\to \D$ is connected if in addition it satisfies
\begin{description}
\item[(3)] for any $d\in \D$ we have $T\downarrow d\neq\emptyset$.
\end{description}
\end{definition}
By abuse of notation, if $T:\C\to\D$ is cellular, we will write $\D-\C$ instead of $\D-T(\C)$.

\begin{remark}\label{cell1}
Suppose that we are given two cellular functors (resp. connected cellular functors) $U:\B\to\C$ and $T:\C\to\D$. Then we claim that the composition $T\circ U$ is also cellular (resp. connected cellular). 
\begin{enumerate}
\item
We first notice that for $b,b'\in\B$
$$
Hom_{\D}(TU(b),TU(b'))\simeq Hom_{\C}(U(b),U(b'))\simeq Hom_{\B}(b,b')
$$ 
so $TU$ is fully faithful. 

\item
Secondly, let $b''\in TU(\B)$. Since $TU(\B)\subset T(\C)$, if $d\in\D-T(\C)$ then $Hom_{\D}(d,b'')=\emptyset$ since $T$ is cellular, and if $d\in T(\C)$, then $Hom_{\D}(d,b'')=\emptyset$ since $T$ is fully faithful and $U$ is cellular. This proves that $TU$ is cellular.

\item
Finally, if both $U$ and $T$ are connected cellular, the categories $T\downarrow d$ and $U\downarrow c$ are not empty for any elements $d\in\D-\C$ and $c\in\C-\B$. So if $d'\in\D-\C$, the category $TU\downarrow d'\neq\emptyset$ and if $d'=T(c')\in T(\C)-T(\B)$, the category $TU\downarrow T(c')\neq\emptyset$. This further proves that $TU$ is connected.
\end{enumerate}
\end{remark}

\begin{definition}[Grothendieck Construction, \cite{Gro}, VI, Sections 8, 9]
Let\ \  $\F: \D\to\mathbf{Cat}$ be a functor. The Grothendieck construction $\D\int\F$ is defined as a category such that:
\begin{enumerate}
\item Objects in $\D\int\F$ are pairs $(d,x)$ such that $d\in\D$ and $x\in\F(d)$;
\item For two objects $(d,x),\ (d',x')\in\D\int\F$, a morphism between them is given by the pair $d\xrightarrow{f}d'$ and $\F(f)(x)\xrightarrow{g}x'$;
\item For three objects $(d,x),\ (d',x')$ and $(d'',x'')$, the composition of morphisms is given by $d\xrightarrow{f}d'\xrightarrow{f'}d''$ and\  $\F(f'\circ f)(x)\xrightarrow{\F(f')(g)}\F(f')(x')\xrightarrow{g'}x''$. 
\end{enumerate}
\end{definition}
\

\begin{remark}\label{rem1.3}
One notices that a morphism $(f,g)\in Mor(\D\int\F)$ has a canonical factorization
$$(f,g)=(1, g)\circ(f, 1_{\F(f)(x)})$$
where $(f, 1_{\F(f)(x)})(d,x)=(d', \F(f)(x))$ and $(1_{d'},g)(d', \F(f)(x))=(d', x')$. In other words, a morphism  can be expressed by the commutative triangle:
$$
\xymatrix{
(d,x)\ar@{|->}[rr]^{(f,g)}\ar@{|->}[rd]_{(f, 1_{\F(f)(x)})} & & (d',x')\\
 & (d',\F(f)(x))\ar@{|->}[ru]_{(1,g)}
}
$$
\end{remark}
\

Let $T:\C\to\D$ be a functor between two small categories. Recall that the functor $\F_T: \D\to\mathbf{Cat}$ sends $d\in\D$ to $T\downarrow d$.
The Grothendieck construction $\D\int\F_T$ has objects $(d,x)$ such that $d\in\D,\ x\in\F_T(d)=T\downarrow d$. A morphism $(d,x)\to(d',x')\in\D\int\F_T$ is given by a morphism $f: d\to d'\in\D$ and $g: \F_T(f)(x)\to x'\in\F_T(d')=T\downarrow d'$, i.e, a commutative diagram:
$$
\xymatrix{
T(c)\ar[r]\ar[d]_{x} & T(c')\ar[d]^{x'}\\
d\ar[r]_f & d'
}
$$

\begin{definition}
We say that a (naturally) commutative diagram of categories and functors is homotopy cocartesian if it is so after applying the nerve functor (cf. \cite[definition 2.3.5]{K}).
\end{definition}

\begin{proposition}[\cite{K}, proposition 2.3.5]\label{cocart}
If $T:\C\to\D$ is cellular, the naturally commutative diagram of categories
$$
\xymatrix{
(\D-\C)\int\F_T\ar[r]^{\ \ \ \ \ p}\ar[d]_{\epsilon} & \C\ar[d]^T\\
\D-\C\ar[r]_{\iota} & \D
}
$$
is homotopy cocartesian. Here, $\iota$ is inclusion, $p$ is the projection to $\C$ and $\epsilon$ is the augmentation induced by $\F_T\to \star$. 
\end{proposition}

\begin{remark}
We notice that if $\D=\C\coprod \C'$ is a category of the disjoint union of two categories $\C$ and $\C'$, then the diagram of proposition \ref{cocart} will be trivial since $\forall\ c'\in\C',\ \F_T(c')=\emptyset$. For a connected cellular functor, this trivial situation is avoided. We also remark that in Quillen's case when $\D=Q_n$, $\C=Q_{n-1}$ and $\C\xrightarrow{T}\D$ is the natural inclusion functor, we always have $\F_T(d)\neq\emptyset$ for any $d\in\D-\C$. 
\end{remark}
\

\subsection{Thomason's Theorem}
Recall that whenever we are given a functor $\F:\D\to\mathbf{Cat}$, we can define the singular chain complex $C_*(\D,\F)$ (c.f, definition \ref{singchain}).
\begin{definition}\label{defiber}
Let $\F: \D\to\mathbf{Cat}$ be a functor such that $\F(d)\neq\emptyset$ for all $d\in\D$, then the relative reduced chain complex $C_*(\D, \widetilde{\F})$ is defined as:
$$C_*(\D, \widetilde{\F}):=Ker(C_*(\D, \F)\to C_*(\D)).$$
If $\F(d)=\emptyset$ for some $d\in\D$, then $C_*(\D, \widetilde{\F})$ is defined to be the homotopy fiber of $C_*(\D, \F)\to C_*(\D)$ in the derived category of Abelian groups $\mathbf{D}(Ab)$.
\end{definition}

\begin{theorem}[Thomason's Theorem, \cite{T}, theorem 1-2]\label{thomason}
Let $\C$ be a category and $\F: \C\to\mathbf{Cat}$ be a functor. There is a weak equivalence:
$$\delta N(\C, \F)\xrightarrow{\sim} N\left(\C\int\F\right),$$
sending 
$$(c_0\xrightarrow{f_1}\cdots\xrightarrow{f_n}c_n,\  x_0\xrightarrow{g_1}\cdots\xrightarrow{g_n}x_n); \ \ \ c_i\in\C,\ x_i\in\F(c_0)$$
to
$$
\left((c_0,y_0=x_0)\xrightarrow{h_1}\cdots\xrightarrow{h_n}(c_n,y_n)\right)
$$
with $y_i=\F(f_i\cdots f_1)x_i\in\F(c_i)$ for $0<i\leq n$ and $h_i=(f_i,\F(f_i\cdots f_1)g_i)$.
\end{theorem}
\

According to Thomason's theorem, for any small category $\D$ and any functor $\F:\D\to\mathbf{Cat}$, the complex $C_*(\D,\F)$ is quasi-isomorphic to the complex $C_*(\D\int\F)$, so we have:
\begin{theorem}[\cite{K}, theorem 2.3.6]\label{homotopycone}
If $T:\C\to\D$ is a cellular functor, there exists an distinguished triangle in the derived category of abelian groups:
$$
C_*(\D-\C, \widetilde{\F_T})\xrightarrow{\widetilde{p}_*} C_*(\C)\xrightarrow{T_*} C_*(\D)\xrightarrow{\partial_T} C_*(\D-\C, \widetilde{\F_T})[1].
$$
\end{theorem}
This theorem gives us a description of the homotopy mapping cone of $T_*$ as $C_*(\D-\C, \widetilde{\F_T})[1]$.
\\

\subsection{Induced Cellular Functor}

Now let $T:\C\to\D$ be a cellular functor, $\iota:\D-\C\to\D$ is the natural inclusion and $R: \D\to\mathbf{Cat}$ be an arbitrary functor. We get the induced functors $R\circ T: \C\to\mathbf{Cat}$ and $R\circ\iota: \D-\C\to\mathbf{Cat}$. So there are Grothendieck constructions $\D\int R$, $\C\int R\circ T$ and $(\D-\C)\int R\circ\iota$.

\begin{lemma}\label{inducedcell}
The functor
$$T\int R:\  \C\int R\circ T\to \D\int R;\ \ \ \ \ (c,x)\mapsto (T(c), x)$$
sending $(c\xrightarrow{f} c', RT(f)(x)\to x')\in Mor(\C\int RT)$ to $(T(c)\xrightarrow{T(f)}T(c'), RT(f)(x)\to x')\in Mor(\D\int R)$ is cellular. 
\end{lemma}
\begin{proof}: Since $T$ is fully faithful, $T\int R$ is also fully faithful by definition. The category $\D\int R-(T\int R)(\C\int RT)$ can be written as $(\D-\C)\int R\circ\iota$ where $\iota: \D-\C\to \D$. So for any $(d,x)\in(\D-\C\int R\circ\iota)$ and $(c,x')\in\C\int RT$, we have $(\D\int R)((d,x), (c,x'))=\emptyset$, because $\D(d,c)=\emptyset$. So $T\int  R$ is cellular.
\end{proof}

\begin{remark}
We notice that even if $T:\C\to\D$ is a connected cellular functor, the induced functor $T\int R: \C\int R\circ T\to \D\int R$ is not necessarily connected. Because for an arbitrary object $(d,y)\in\D\int R$, we have:
$$
\left( T\int R \right)\downarrow (d,y)=\{(c,x, f, \phi)\ |\ (c,x)\in\C\int R\circ T,\   T(c)\xrightarrow{f} d,\  Rf(x)\xrightarrow{\phi} y\}.
$$
So we see that even if there always exists $f$, the existence of $\phi\in R(d)$ is not guaranteed.
\end{remark}
\

\subsection{First Reduction of the Problem}
Recall that (c.f, \cite[Section 2]{K}) the rank spectral sequence is constructed under the background of any triangulated category $\mathcal{T}$ with countable direct sums as follows. We take 
$$
C_0\xrightarrow{i_1}\cdots\xrightarrow{i_n}C_n\xrightarrow{i_{n+1}}\cdots
$$
being a sequence of objects in $\mathcal{T}$ and $C$ as a homotopy colimit of the $C_n$ (please refer to \cite{BN} for its definition). For each $n$, we choose a mapping cone $C_{n/n-1}$ of $i_n$ so that there is an distinguished triangle
$$
C_{n-1}\xrightarrow{i_n}C_n\xrightarrow{j_n}C_{n/n-1}\xrightarrow{k_n}C_{n-1}[1].
$$
Let $H_*: \mathcal{T}\to\mathcal{A}$ be a (co)homological functor taking values in an abelian category such that $H_*$ commutes with countable direct sums. Then the exact couple
$$
D_{p,q}:=H_{p+q}(C_p),\ \ \ \ \ \ \ E_{p,q}:=H_{p+q}(C_{p/p-1})
$$
gives rise to a spectral sequence $E_{p,q}=E^1_{p,q}\Rightarrow H_{p+q}(C)$ where $C$ is the homotopy colimit of the $C_n$.
Let $C_{p/p-2}$ be a cone of $i_pi_{p-1}: C_{p-2}\to C_p$, then we get a commutative diagram of distinguished triangles
\begin{equation}\label{tri1}
\xymatrix{
C_{p-2}\ar[d]_{i_{p-1}}\ar@{=}[r] & C_{p-2}\ar[d]^{i_pi_{p-1}} & & \\
C_{p-1}\ar[r]^{i_p}\ar[d]_{j_{p-1}} & C_p\ar[r]^{j_p}\ar[d] & C_{p/p-1}\ar[r]^{k_p}\ar[d]^= & C_{p-1}[1]\ar[d]^{j_{p-1}[1]}\\
C_{p-1/p-2}\ar[r]_{\bar{i}_p} & C_{p/p-2}\ar[r]_{\bar{j}_p} & C_{p/p-1}\ar[r]_{\bar{k}_p} & C_{p-1/p-2}[1]
}
\end{equation}
The differential $d^1_{p,q}$ is the boundary map $\bar{k}_{p,n}$ with $n=p+q$ of the long distinguished sequence 
$$
H_{p+q}(C_{p-1/p-2})\xrightarrow{\bar{i}_{p,n}}H_{p+q}(C_{p/p-2})\xrightarrow{\bar{j}_{p,n}}H_{p+q}(C_{p/p-1})\xrightarrow{\bar{k}_{p,n}}H_{p+q-1}(C_{p-1/p-2})
$$
associated with the bottom distinguished triangle of the above diagram. We notice that $d^1$ is unique up to isomorphisms between the chosen mapping cones $C_{n/n-1}$, i.e, if we choose another mapping cone $C'_{n/n-1}$ and an isomorphism $C'_{n/n-1}\to C_{n/n-1}$ fitting into the following commutative diagram of distinguished triangles
$$
\xymatrix{
C_{n-1}\ar[r]^{i_n}\ar@{=}[d] & C_n\ar[r]^{j_n}\ar@{=}[d] & C_{n/n-1}\ar[r]^{k_n} & C_{n-1}[1]\ar@{=}[d] \\
C_{n-1}\ar[r]_{i_n} & C_n\ar[r]_{j'_n} & C'_{n/n-1}\ar[r]_{k'_n}\ar[u]^{\simeq} & C_{n-1}[1]
}
$$
then the different choice of mapping cone gives us another map $\bar{k}'_p: C'_{p/p-1}\to C'_{p-1/p-2}[1]$
and a commutative diagram determined by this isomorphism of mapping cones. It follows that up to the isomorphisms of mapping cones, $d^1$ is unique.
\\

By proposition \ref{cocart}, for a cellular functor $T: \C\to\D$, we choose the mapping cone of $C_*(\C)\to C_*(\D)$ as $C_*(\D-\C, \widetilde{\F_T})[1]$.
Consider two composable cellular functors $\B\xrightarrow{U}\C\xrightarrow{T}\D$. By remark \ref{cell1}, the composition $T\circ U$ is still cellular. Replacing the abstract objects by chain complexes in diagram (\ref{tri1}), we get the following commutative diagram:
\begin{equation}\label{diag1.4}
\xymatrix{
C_*(\B)\ar[d]\ar@{=}[r] & C_*(\B)\ar[d] & & \\
C_*(\C)\ar[r]\ar[d] & C_*(\D)\ar[r]\ar[d] & C_*(\D-\C, \widetilde{\F_T})[1]\ar[r]\ar[d]^= & C_*(\C)[1]\ar[d]\\
C_*(\C-\B,\widetilde{\F_U})[1]\ar[r] & C_*(\D-\B,\widetilde{\F_{TU}})[1]\ar[r] & C_*(\D-\C,\widetilde{\F_T})[1]\ar[r]_{d^1[1]} & C_*(\C-\B,\widetilde{\F_U})[2]
}
\end{equation}
We see that $d^1$ is the composition of
$$
C_*(\D-\C,\widetilde{\F_T})\to C_*(\C)\to C_*(\C-\B,\widetilde{\F_U})[1].
$$
In the following sections, we will put this composition into suitable distinguished triangles so that it can be explicitly calculated.
\\

\section{Associativity of the Grothendieck Construction}
In this section, we will give an explicit description of the double Grothendieck construction and prove the associativity of Grothendieck constructions.

\subsection{The category $(\D\int\F)\int\G$}
For a functor $\F:\D\to\mathbf{Cat}$, we have already had a Grothendieck construction $\D\int \F$. 

\begin{lemma}
A functor $\G: \D\int\F\to\mathbf{Cat}$ is equivalent to the following data: 
\begin{enumerate}
\item For any object $(d,x)\in Obj(\D\int\F)$,\ \  a small category $\G(d,x)$,
\item For any $f: d\to d'\in Mor(\D)$, a functor $\G(f, 1_{\F(f)(x)}): \G(d,x)\to \G(d', \F(f)(x))$,
\item For any $g: x\to x'\in Mor(\F(d))$,  a functor $\G(1_d, g): \G(d,x)\to\G(d,x')$.
\end{enumerate}
which verifies the following axioms:
\begin{description}
\item[(i)] $\G(d, 1_x)=Id_{\G(d,x)}$ and $\G(1_d, x)=Id_{\G(d,x)}$.
\item[(ii)] The functors $\G(f, 1_{\F(f)(x)})$ and $\G(1_d, g)$ are natural in $f$ and $g$ respectively. 
\item[(iii)] For any $d\xrightarrow{f}d'\in\D$ and any $x\xrightarrow{g} x'\in\F(d)$ we have a commutative square of categories:
\begin{equation}\label{diag2.1}
\xymatrix{
\G(d,x)\ar[rrr]^{\G(f,1_{\F(f)(x)})}\ar[d]_{\G(1_d,g)} & & & \G(d', \F(f)(x))\ar[d]^{\G(1_{d'}, \F(f)(g))}\\
\G(d,x')\ar[rrr]_{\G(f, 1_{\F(f)(x')})} & & & \G(d', \F(f)(x'))
}
\end{equation}
\end{description}

\end{lemma}

\begin{proof}:
If we have a functor $\G: \D\int\F\to\mathbf{Cat}$, it clearly verifies the condition (1) and axiom (i). By remark \ref{rem1.3}, any morphism $(f,g)=(d\xrightarrow{f}d', \F(f)(x)\xrightarrow{g}x'): (d,x)\to (d',x')\in Mor(\D\int\ F)$ can be factorized as: $(f,g)=(1_d,g)\circ(f, 1_{\F(f)(x)})$, so we have a facorization
$$\G(f,g):=\G(1_{d'},g)\circ\G(f, 1_{\F(f)(x)}): \G(d,x)\to \G(d',x')\in\mathbf{Cat},$$
so we see that $\G$ verifies condition (2), (3) and the functoriality axiom (ii). 
Moreover, for any $f: d\to d'\in\D$ and $g: x\to x'\in\F(d)$, a morphism $(d,x)\to (d', \F(f)(x'))\in\D\int\F$ can be factorized as $(1_{d'}, \F(f)(g))\circ (f, 1_{\F(f)(x)})=(1_d,g)\circ(f, 1_{\F(f)(x')})$. Applying $\G$ to this factorization, we find $\G$ verifies axiom (iii). 
\\

On the other hand, suppose we are given conditions (1)-(3) that verify the three axioms, we will verify that $\G$ is a well-defined functor. The condition (1) means we have done on object level, so it leaves us to verify on the morphism level, more precisely, composition of morphisms.
\\

Note that the functoriality of $\G(f, 1_{\F(f)(x)})$  can be shown by a commutative diagram:
\begin{equation}\label{diag2.2}
\xymatrix{
\G(d,x)\ar[rr]^{\G(f, 1_{\F(f)(x)})}\ar[rrdd]_{\G(f'\circ f, 1_{\F(f'\circ f)(x)})} & & \G(d', \F(f)(x))\ar[dd]^{\G(f', 1_{\F(f'\circ f)(x)})}\\
& & \\
& & \G(d'', \F(f')\circ\F(f)(x))
}
\end{equation}
such that $d\xrightarrow{f}d'\xrightarrow{f'}d''\in Mor(\D)$.
For any $(d,x)\xrightarrow{(f,g)}(d',x')\in\D\int\F$, there exists a functor
$$\G(f,g):=\G(1_d,g)\circ\G(f, 1_{\F(f)(x)}): \G(d,x)\to \G(d',x')\in\mathbf{Cat},$$
that is, for any $x\in\F(d)$ and $x'\in\F(d')$, we have a commutative diagram of functors:
$$
\xymatrix{
\G(d,x)\ar[rrr]^{\G(f, 1_{\F(f)(x)})}\ar[rrrd]_{\G(f,g)} & & & \G(d', \F(f)(x))\ar[d]^{\G(1_{d'},g)}\\
 & & & \G(d',x')
}
$$
So, combined with diagram (\ref{diag2.1}) and (\ref{diag2.2}), compositions can be described as:
$$
\xymatrix{
\G(d,x)\ar[rrr]^{\G(f, 1_{\F(f)(x)})}\ar[rrrd]_{\G(f,g)} & & & \G(d', \F(f)(x))\ar[d]^{\G(1_{d'},g)}\ar[rrr]^{\G(f', 1_{\F(f'\circ f)(x)})} & & & \G(d'', \F(f'\circ f)(x))\ar[d]^{\G(1_{d''}, \F(f')(g))}\\
 & & & \G(d',x')\ar[rrr]^{\G(f', 1_{\F(f')(x')})}\ar[rrrd]_{\G(f',g')} & & & \G(d'',\F(f')(x'))\ar[d]^{\G(1_{d''},g')}\\
 & & & & & & \G(d'',x'')
}
$$
which means 
$$\G(1_{d''}, g'\circ \F(f')(g))\circ\G(f'\circ f, 1_{\F(f'\circ f)(x)})=\G(1_{d''},g')\circ\G(f', 1_{\F(f')(x')})\circ \G(1_{d'},g)\circ\G(f, 1_{\F(f)(x)}).$$
Thus, by remark \ref{rem1.3}, we have:
$$\G(f'\circ f, g'\circ \F(f')(g))=\G(f',g')\circ\G(f,g).$$ 
So $\G$ is a well-defined functor. 

\end{proof}
\

We have a Grothendieck construction $(\D\int\F)\int\G$ such that the objects are $((d,x), \alpha)$ with $(d,x)\in \D\int\F$ and $\alpha\in\G(d,x)$. For two objects $((d,x),\alpha)$ and $((d',x'), \alpha')$, a morphism between them is defined as $(d\xrightarrow{f}d', \F(f)(x)\xrightarrow{g}x')\in\D\int\F$ and $\G(f,g)(\alpha)\to\alpha'$. For three objects $((d,x), \alpha),\ ((d',x'), \alpha')$ and $((d'',x''), \alpha'')$, the composition of two morphisms is the data
$$d\xrightarrow{f}d'\xrightarrow{f'}d'';\ \ \ \ \F(f'f)(x)\xrightarrow{\F(f')(g)} \F(f')(x')\xrightarrow{g'} x''$$
and
$$\G(f',g')\circ\G(f,g)(\alpha)\to\G(f',g')(\alpha')\to\alpha''.$$
\\

\subsection{The category $\D\int(\F\int\G)$}
Let $\D$ be a category, $\F:\D\to\mathbf{Cat}$ be a functor and $\G: \D\int\F\to\mathbf{Cat}$ be another functor. We will first construct a functor $\F\int\G: \D\to\mathbf{Cat}$. We define (c.f, diagram of remark \ref{rem1.3}) 
$$
\G_f(x):=\G(f, 1_{\F(f)(x)}),\ \ \ \ \ \ \ d\xrightarrow{f}d'\in\D,\ \ x\in \F(d)
$$
and $\G_d:=\G(d,-): \F(d)\to\mathbf{Cat}$ sends $x\mapsto\G(d,x)$.
\\

Let us now define a functor $\F\int\G: \D\to\mathbf{Cat}$:

\begin{description}

\item[On objects] $d\mapsto \F(d)\int\G_d=\{(x,\alpha)\ |\ x\in\F(d), \alpha\in\G_d(x)\}$. An object $(x,\alpha)\in\F(d)\int\G_d$ is given by $x\in\F(d)$ and $\alpha\in\G_d(x)$. For two objects $(x,\alpha),\ (x',\alpha')$ in the category $\F(d)\int\G_d$, a morphism between them is denoted by $(h,l)$ which is given by $x\xrightarrow{h}x'\in\F(d)$ and the morphism $\G_d(h)(\alpha)\xrightarrow{l} \alpha'\in \G_d(x')$.
\\

\item[On morphisms] If $d\xrightarrow{f}d'\in\D$ is a morphism, then $(\F\int\G)(f)$ is the functor defined by:
$$(\F\int\G)(f): \F(d)\int\G_d\to\F(d')\int\G_{d'};\ \ \ \ \ \ (x, \alpha)\mapsto (\F(f)(x), \G_f(x)(\alpha)).$$
Moreover, if $(h,l): (x,\alpha)\to (x', \alpha')$ is a morphism in $\F(d)\int\G_d$, we define $(\F\int\G)(f)(h,l)=(\F(f)(h), \G_f(x')(l))$ as 
$$
\F(f)(h): \F(f)(x)\to\F(f)(x');\ \ \ \ \ \ \G_f(x')(l): \G_{f}(x')\circ\G_d(h)(\alpha)\to \G_f(x')(\alpha').
$$
If $(x,\alpha)\xrightarrow{(h,l)} (x',\alpha')\xrightarrow{(h',l')} (x'',\alpha'')$ is a chain of morphisms in $\F(d)\int\G_d$, then we have the compositions: $\F(f)(x)\to \F(f)(x')\to \F(f)(x'')$
and
$$\G_f(x'')\circ\G_d(h')\circ\G_d(h)(\alpha)\xrightarrow{\G_f(x'')\circ\G_d(h')(l)}\G_f(x'')\circ \G_{d}(h')(\alpha')\xrightarrow{\G_f(x'')(l')}\G_f(x'')(\alpha''),$$
which equals $\G_f(x'')(l'\circ l)$.
So we have 
$$(\F\int\G)(f)((h',l')\circ (h,l))=(\F(f)(h'\circ h), \G_f(x'')(l'\circ l)).$$
Moreover, if $1=(1_x, 1_{\alpha})$, then $(\F\int\G)(f)(1)=id$. Therefore, $(\F\int\G)(f)$ is a morphism from $\F(d)\int\G_{d'}\to \F(d')\int\G_{d'}$ in $\mathbf{Cat}$ for any $f:d\to d'\in\D$.
\\

\item[Compositions] 
If $d\xrightarrow{f} d'\xrightarrow{f'}d''$ are two composable morphisms in $\D$, we have:
$$
(\F\int\G)(f'\circ f): \F(d)\int\G_d\to \F(d'')\int\G_{d''}
$$
which sends $(x,\alpha)\in \F(d)\int\G_d$ to $(\F(f'\circ f)(x), \G_{f'\circ f}(x)(\alpha))$ which equals to
$$
(\F(f')\circ\F(f)(x),\G_{f'}(\F(f)(x))\circ\G_f(x)(\alpha))=(\F\int\G(f'))\circ(\F\int\G(f))(x,\alpha).
$$

\item[Identity] For any identity morphism $1\in\D$, we have
$$(\F\int\G)(1)(x,\alpha)=(\F(1)(x), \G_1(x)(\alpha))=(x, \alpha),$$
which means $(\F\int\G)(1)=Id$.
\end{description}
\

Thus we see that $\F\int\G$ is a well-defined functor, so there exists a Grothendieck construction $\D\int(\F\int\G)$.
\\

\subsection{Identification of the two Constructions}

\begin{proposition}\label{associate}
The Grothendieck construction is associative, that is to say, under our assumption, we have an isomorphism of categories:
$$(\D\int\F)\int\G=\D\int(\F\int\G).$$
\end{proposition}
\begin{proof}: 
The left hand side has objects 
$$((d,x), \alpha),\ \ \ (d,x)\in\D\int\F,\ \alpha\in\G(d,x),$$
and the right hand side has objects:
$$(d, (x,\alpha)),\ \ \ d\in\D,\ x\in\F(d),\ \alpha\in\G(d,x).$$
So the two constructions coincide on object level.
\\

On morphism level, the morphisms of $(\D\int\F)\int\G$ are described above. On the other hand, a morphism $(d, (x,\alpha))\to (d', (x',\alpha'))$ in $\D\int(\F\int\G)$ consists of:
$$d\xrightarrow{f} d',\ \ \ \ \  (\F\int\G)(f)(x,\alpha)\to (x',\alpha'),$$
where $(\F\int\G)(f)(x,\alpha)=(\F(f)(x), \G(f, 1_{\F(f)(x)})(\alpha))\in\G(d',x')$. Thus $(\F\int\G)(f)(x,\alpha)\to (x',\alpha')$ consists of two morphisms
$$\F(f)(x)\xrightarrow{g}x';\ \ \ \ \ \ \G(f,g)(\alpha)=\G(1,g)\circ\G(f, 1_{\F(f)(x)})(\alpha)\to\alpha'.$$
So we can see that the two constructions coincide on morphism level as well. It is easy to verify that the compositions of two sides match, so these two categories are isomorphic. 

\end{proof}
\

\section{Further Reduction using Grothendieck Constructions}

\subsection{Background}
We consider the composable connected cellular functors 
\begin{equation}\label{I}
\B\xrightarrow{U}\C\xrightarrow{T}\D.
\end{equation}

Recall that by the construction of the rank spectral sequence (Chapter 1.4) we have short exact sequences
\begin{equation}\label{Ia}
0\to C_*(\B)\xrightarrow{U} C_*(\C)\to C_*(\C)/C_*(\B)\to 0
\end{equation}
\begin{equation}\label{Ib}
0\to C_*(\C)\xrightarrow{T} C_*(\D)\to C_*(\D)/C_*(\C)\to 0,
\end{equation}
which induce morphisms in the derived category of abelian groups $\mathbf{Ab}$:
\begin{equation}\label{Ic}
C_*(\D)/C_*(\C)\to C_*(\C)[1]\to C_*(\C)/C_*(\B)[1].
\end{equation}
By theorem \ref{homotopycone}, we choose the following mapping cones
$$
C_*(\D)/C_*(\C)=C_*(\D-\C,\widetilde{\F_T})[1],\ \ \ \ \ \ \ C_*(\C)/C_*(\B)=C_*(\C-\B,\widetilde{\F_U})[1].
$$
Together with diagram (\ref{diag1.4}), we see that in order to calculate $d^1$ of the rank spectral sequence, we calculate the composition of (\ref{Ic}).
\\

In this section, we put (via quasi-isomorphisms) the composition of (\ref{Ic}) into certain distinguished triangles that come from distinguished triangles of functors (coefficients) defined over different categories: firstly $\D$, then $\D-\B$. 
Notice that every short exact sequence is naturally a distinguished triangle, so we say that two short exact sequences are quasi-isomorphic if corresponding distinguished triangles are. In this case, we say that one short exact sequence can be replaced (up to quasi-isomorphism) by the other.
It turns out that $d^1$ of rank spectral sequence can be calculated via the distinguished triangles of functors we constructed.
\\

\subsection{Pass to $\D$-chains}
We define a functor 
$$
U_*: \B\to \D\int\F_T;\ \ \ \ \ b\mapsto (b=b).
$$
Then there exists a Grothendieck construction $\D\int\F_T\int\F_{U_*}$ whose objects can be written as $b\to c\to d$ with $b\in\B$ and $c\to d\in\D\int\F_T$. Notice that for a fixed $c\to d\in\F_T(d)$, the categories $\F_{U_*}(c\to d)$ and $\F_U(c)$ are canonically isomorphic.
\begin{lemma}\label{lemma1}
The categories $\D\int\F_T\int\F_{U_*}$ and $\D\int\F_{TU}$ are homotopy equivalent.
\end{lemma}
\begin{proof}
We define two functors
$$
p: \D\int\F_T\int\F_{U_*}\to \D\int\F_{TU};\ \ \ \ \ (b\to c\to d)\mapsto(b\to d)
$$
and
$$
s: \D\int\F_{TU}\to \D\int\F_T\int\F_{U_*};\ \ \ \ \ (b\to d)\mapsto(b=b\to d)
$$
Since $s$ is left adjoint to $p$, they induce homotopy equivalences between two corresponding categories.

\end{proof}
\

We notice that this lemma actually says that $p$ and $s$ make the functors $\F_T\int\F_{U_*}$ and $\F_{TU}$ naturally homotopy equivalent over the category $\D$, i.e, after applying to $\D$ the resulting categories are homotopy equivalent.
\\

It is easy to see that the projection functor $\D\int\F_T\to\C$ (resp. $\D\int\F_{TU}\to\B$) admits a left adjoint $c\to [c=c]$ (resp. $b\to[b=b]$), so there exist a pair of homotopy equivalences between the categories $\D\int\F_T$ and $\C$ (between $\D\int\F_{TU}$ and $\B$).
Together with Thomason's theorem (theorem \ref{thomason}), we obtain canonical homotopy equivalences
\begin{equation}\label{II}
N(\D, \F_T)\approx N\C,\ \ \ \ \ \ \ \ \ \ N(\D, \F_{TU})\approx N\B.
\end{equation}
Then, up to canonical homotopy equivalences, (\ref{I}) becomes
\begin{equation}\label{II'}
\D\int\F_{TU}\to \D\int\F_T\to \D,
\end{equation}
and (\ref{Ib}) can be replaced by (notice that $T$ is connected cellular)
\begin{equation}\label{IIb}
0\to C_*(\D, \widetilde{\F_T})\to C_*(\D, \F_T)\to C_*(\D)\to 0.
\end{equation}
\begin{definition}[Compare to \cite{K}, 2.2.2]
Suppose that $U:\B\to\C$ is connected cellular. We define the relative reduced chain complex
$$
C_*(\D, \F_T\int\widetilde{\F_{U_*}})=Ker( C_*(\D, \F_T\int\F_{U_*})\to C_*(\D, \F_T)).
$$
\end{definition}
Therefore, use lemma \ref{lemma1}, we see that (\ref{Ia}) and (\ref{Ic}) can be replaced by
\begin{equation}\label{IIa}
0\to C_*(\D, \F_T\int\widetilde{\F_{U_*}})\to C_*(\D, \F_T\int\F_{U_*})\to C_*(\D, \F_T)\to 0.
\end{equation}
and
\begin{equation}\label{IIc}
C_*(\D, \widetilde{\F_T})\to C_*(\D, \F_T)\to C_*(\D,\F_T\int\widetilde{\F_{U_*}})[1]\xrightarrow{+1}
\end{equation}
In particular, the above triangle is induced by the distinguished triangle of their coefficients:
\begin{equation}\label{IIc'}
C_*(\widetilde{\F_T})\to C_*(\F_T)\to C_*(\F_T\int\widetilde{\F_{U_*}})[1]\xrightarrow{+1}
\end{equation}
\

\subsection{distinguished Triangles of Functors/Coefficients}
There exists a commutative diagram of short exact sequences
\begin{equation}\label{seq2}
\xymatrix{
0\ar[r] & C_*(\D, \F_T\int\widetilde{\F_{U_*}})\ar@{=}[d]\ar[r] & C_*(\D, \widetilde{\F_T\int\F_{U_*}})\ar[r]\ar[d] & C_*(\D, \widetilde{\F_T})\ar[r]\ar[d] & 0\\
0\ar[r] & C_*(\D, \F_T\int\widetilde{\F_{U_*}})\ar[r]\ar[d] & C_*(\D, \F_T\int\F_{U_*})\ar[r]\ar[d] & C_*(\D, \F_T)\ar[r]\ar[d] & 0\\
0\ar@{=}[r] & 0\ar[r] & C_*(\D)\ar@{=}[r] & C_*(\D)\ar[r] & 0
}
\end{equation}
such that all sequences are induced by the short exact sequences of their coefficients.
Thus, the composition of (\ref{IIc}) fits into the distinguished triangle
$$
C_*(\D, \F_T\int\widetilde{\F_{U_*}})\to C_*(\D,\widetilde{\F_T\int\F_{U_*}})\to C_*(\D,\widetilde{\F_T})\xrightarrow{+1}C_*(\D, \F_T\int\widetilde{\F_{U_*}})[1]
$$
and the composition of
(\ref{IIc'}) fits into the distinguished triangle of coefficients
\begin{equation}\label{IIc''}
C_*(\F_T\int\widetilde{\F_{U_*}})\to C_*(\widetilde{\F_T\int\F_{U_*}})\to C_*(\widetilde{\F_T})\xrightarrow{+1}C_*(\F_T\int\widetilde{\F_{U_*}})[1].
\end{equation}
This triangle is distinguished in the sense that it gives a term-wise distinguished triangle, i.e,  after applying to any object $d\in\D$ the resulting triangle is distinguished. For simplicity, we will just say that this sequence is a distinguished triangle of coefficients or functors without mentioning which category it applies. 
\\

\subsection{Pass to $(\D-\B)$-chains}
There exists a commutative diagram of distinguished triangles
\begin{equation}\label{III}
\xymatrix{
C_*(\D-\B, \F_T\int\widetilde{\F_{U_*}})\ar[r]\ar[d] & C_*(\D-\B, \widetilde{\F_T\int\F_{U_*}})\ar[r]\ar[d] & C_*(\D-\B, \widetilde{\F_T})\ar[r]^{\ \ \ \ \ \ +1}\ar[d] & \\
C_*(\D, \F_T\int\widetilde{\F_{U_*}})\ar[r] & C_*(\D, \widetilde{\F_T\int\F_{U_*}})\ar[r] & C_*(\D, \widetilde{\F_T})\ar[r]^{\ \ \ \ \ \ \ +1} & 
}
\end{equation}
where all vertical morphisms are quasi-isomorphisms (we use \cite[Prop. 2.3.4]{K} to the case $\phi: \F\Rightarrow *$ and $TU: \B\to\D$, where we replace $\F$ by corresponding functors appearing in the above diagram).
It follows that $d^1$ can be calculated by the top distinguished triangle of diagram (\ref{III}) which is induced by the distinguished triangle of coefficients (\ref{IIc''}).

\begin{remark}\label{lemma1'}
Use the proof of lemma \ref{lemma1}, we can show that $p$ and $s$ give homotopy equivalences between $(\D-\B)\int\F_T\int\F_{U_*}$ and $(\D-\B)\int\F_{TU}$. So $\F_T\int\F_{U_*}$ and $\F_{TU}$ are naturally homotopy equivalent over the category $\D-\B$.
\end{remark}
\

We define $\C':=\C-\B$, $\D':=\D-\B$ and $T'=T|_{\C'}: \C'\to\D'$. Notice that if $T$ is cellular then so does $T'$.  
\\

\subsection{Replace $T$ by $T'$}
Suppose that $T:\C\to\D$ and $T':\C'\to\D'$ are connected cellular. 
We first notice that for any $d\in\D'$ there exists a functor 
$$
\F_{U_*}|_{\F_{T'}(d)}: \F_{T'}(d)\to \mathbf{Cat}
$$
so that we have a Grothendieck construction $\F_{T'}(d)\int\F_{U_*}$ whose objects are $b\to c\to d$ with $b\in\B$ and $c\in\C-\B$.
\\

Let $d\in\D'$. The inclusion functor $i: \F_{TU}(d)\hookrightarrow\F_T(d)$ is connected cellular and  
$$
\F_T(d)-\F_{TU}(d)=\F_{T'}(d).
$$
Applying \cite[Proposition 2.3.4]{K} to the cellular functor $i$ and the natural transformation $\F_{U_*}\Rightarrow \star$ and use Thomason's formula, we can construct a homotopy cocartesian diagram
\begin{equation}\label{diag2}
\xymatrix{
\F_{T'}(d)\int\F_{U_*}\ar[rr]\ar[d]& & \F_T(d)\int\F_{U_*}\ar[d]\\
\F_{T'}(d)\ar[rr] & & \F_T(d)
}
\end{equation}
which gives a canonical quasi-isomorphism
\begin{equation}\label{eq8}
C_*(\F_{T'}(d)\int\widetilde{\F_{U_*}})\xrightarrow{\sim} C_*(\F_T(d)\int\widetilde{\F_{U_*}}),\ \ \ \ \forall d\in\D'.
\end{equation}
Combined with the top distinguished triangle of diagram (\ref{III}), we obtain an distinguished triangle
\begin{equation}\label{IV}
C_*(\D-\B, \F_{T'}\int\widetilde{\F_{U_*}})\to C_*(\D-\B, \widetilde{\F_T\int\F_{U_*}})\to C_*(\D-\B, \widetilde{\F_T})\xrightarrow{+1}
\end{equation}
which is induced by the distinguished triangle of coefficients/functors
\begin{equation}\label{IIc'''}
C_*(\F_{T'}\int\widetilde{\F_{U_*}})\to C_*(\widetilde{\F_T\int\F_{U_*}})\to C_*(\widetilde{\F_T})\xrightarrow{+1}C_*(\F_{T'}\int\widetilde{\F_{U_*}})[1].
\end{equation}
Moreover, according to remark \ref{lemma1'}, (\ref{IV}) is isomorphic to the distinguished triangle
\begin{equation}\label{IV'}
C_*(\D-\B, \F_{T'}\int\widetilde{\F_{U_*}})\to C_*(\D-\B, \widetilde{\F_{TU}})\to C_*(\D-\B, \widetilde{\F_T})\xrightarrow{+1}
\end{equation}
which is induced by the distinguished triangle of coefficients/functors
\begin{equation}\label{IIc''''}
C_*(\F_{T'}\int\widetilde{\F_{U_*}})\to C_*(\widetilde{\F_{TU}})\to C_*(\widetilde{\F_T})\xrightarrow{+1}C_*(\F_{T'}\int\widetilde{\F_{U_*}})[1].
\end{equation}
\

We consider Quillen's Q-construction: let $A$ be an integral Noetherian domain and $\D=Q_n$ be the full sub-category of $Q^{\tf}(A)$ consisting of torsion-free modules of rank smaller or equal to $n$. Similarly, we define $\B=Q_{n-2}$ and $\C=Q_{n-1}$.
Let $T':\C-\B\hookrightarrow\D-\B$ be the inclusion functor. We notice that for any $d\in \D-\C$ the category $\F_{T'}(d)$ has only objects the admissible monomorphisms $c\rightarrowtail d$ and admissible epimorphisms $d\twoheadrightarrow c$, hence it is discrete. On the other hand, if $d\in\C-\B$ then $\F_{T'}(d)$ is contractible since this category has a terminal object $d=d$. In particular, applying (\ref{IIc''''}) to any $d\in\D-\C$, we get a distinguished triangle
\begin{equation}\label{IIc'''''}
\bigoplus_{c\to d\in\F_{T'}(d)}C_*(\widetilde{\F_U}(c))\to C_*(\widetilde{\F_{TU}}(d))\to C_*(\widetilde{\F_T}(d))
\end{equation}
$$
\xrightarrow{+1}\bigoplus_{c\to d\in\F_{T'}(d)}C_*(\widetilde{\F_U}(c))[1].
$$
Therefore, $d^1$ on coefficients is given by
\begin{equation}\label{coeff}
\widetilde{St}(d)\simeq H_{n-1}(\widetilde{\F_T}(d))\to H_{n-2}(\F_{T'}(d)\int\widetilde{\F_{U_*}})
\end{equation}
$$
=\bigoplus_{c\to d\in\F_{T'}(d)}H_{n-2}(\widetilde{\F_{U_*}}(c\to d))\simeq\bigoplus_{c\to d\in\F_{T'}(d)}\widetilde{St}(c).
$$

\section[Modular Symbols and Tits Buildings]{An Intermezzo of Modular Symbols and Tits Buildings}
We first recall following notations and results given in \cite[Section 2]{KS}. 
We shall work with essentially $4$ categories: 
\begin{itemize}
\item $\Set$, the category of (small) sets.
\item $\Ord$, the category of partially ordered sets. Recall that, as in Quillen \cite{Q1}, we may think of a poset as a category.
\item $\Spl$, the category of abstract simplicial complexes (see definition \ref{defsc} below).
\item $\Top$, the category of topological spaces.
\end{itemize}

There are various functors between these categories: we write
\begin{itemize}
\item $E:\Set \to \Spl$  for the functor which sends a set $X$ to the simplicial complex of nonempty finite subsets of $X$.  
\item $B:\Ord\to \Spl$ for the functor sending a poset to the  simplicial complex of its  totally ordered nonempty finite subsets.
\item $\Simpl:\Spl\to \Ord$ for the functor which associates to a simplicial complex the set of its simplices ordered by inclusion. 
\item $sd: \Spl\to\Spl$ for the first barycentric subdivision functor which is defined as $sd=B\circ\Simpl$.
\item $|\ |:\Spl\to \Top$ for the geometric realisation functor \cite[3.1]{Sp}.
\end{itemize}
\

There exists a homotopy equivalence (c.f, \cite[(3)]{KS})
\begin{equation}\label{e3}
\epsilon_\Gamma:| sd \Gamma|\xrightarrow{\approx} |\Gamma|
\end{equation}
which is natural in $\Gamma\in \Spl$.
For any $S\in \Ord$, $|B(S)|$ is naturally homeomorphic to $|N(S)|$, where $N(S)$ is the nerve of the category $S$; conversely, if $\Gamma\in \Spl$, the relation $B\circ \Simpl=sd$\ and (\ref{e3}) yield a natural  homotopy equivalence $|N(\Simpl(\Gamma))|\xrightarrow{\approx} |\Gamma|$   (compare \cite[p. 89]{Q1}). Thus we can work equivalently with simplicial complexes or posets, and use Quillen's techniques from \cite{Q1} when dealing with the latter. Following the practice in \cite{Q1} and \cite{Q2}, we shall say that a poset, a simplicial complex, or a morphism in $\Ord$ or $\Spl$ have a certain homotopical property if their topological realisations have.
\\

\subsection{(Universal) Modular Symbols}

\begin{definition}[\cite{Sp}, Chapter III, section 1]\label{defsc}
A simplicial complex $K$ consists of a set $\mathcal{V}$ of vertices and a set $\mathcal{S}$ of finite non-empty subsets of $\mathcal{V}$ called simplices such that
\begin{enumerate}
\item Any set consisting of exactly one vertex is a simplex.
\item Any non-empty subset of a simplex is a simplex.
\end{enumerate}
\end{definition}
\

Recall that $\Sigma T(V)$ denotes Quillen's model of the suspension of the Tits building $T(V)$. Let $[n]=\{0,1,\cdots,n\}$ be the ordered set for $n\geq 0$, we denote $\Delta_n$ the standard $n$-simplex defined as $\Delta_n:=E([n])$. 
\begin{definition}[Modular Symbols, c.f, \cite{AR}]\label{ms}
Let $g_0,\cdots,g_{n-1}\in V-\{0\}$ be a collection of non-zero vectors and $sd\Delta_{n-1}$ be the barycentric subdivision of $\Delta_{n-1}$. The vertices of $\partial sd\Delta_{n-1}$ are the proper non-empty subsets of $[n-1]=\{0,1,\cdots,n-1\}$. We denote by $\phi: \partial sd\Delta_{n-1}\to T(V)$ the simplicial map sending each vertex $v\in\partial sd\Delta_{n-1}$ to the subspace of $V$ generated by $\{g_i\ |\ i\in v\}$. Let $\zeta_{n-1}$ be the fundamental class of $H_{n-2}(\partial sd\Delta_{n-1})$, then we call $[g_0,\cdots,g_{n-1}]:=\phi_*(\zeta_{n-1})\in H_{n-2}(T(V))$ a (universal) modular symbol.
\end{definition}
\

\begin{proposition}[c.f, Prop. 2.2, \cite{AR}]\label{propms}
Let $[G]=[g_0,\dots,g_{n-1}]$ be a modular symbol, then:
\begin{enumerate}
\item $[g_0,\dots,g_i,\dots,g_j,\dots,g_{n-1}]=-[g_0,\dots,g_j,\dots,g_i,\dots,g_{n-1}]$ for any $0\leq i\neq j\leq n-1$;
\item $[ag_0,\dots,g_{n-1}]=[g_0,\dots, g_{n-1}]$ for any $a\in K-\{0\}$;
\item if $det(G)=0$ then $[g_0,\dots,g_{n-1}]=0$;
\item if $g_0,\dots,g_n$ be $n+1$ non-zero vectors, then:
             $$
             \sum_{i=0}^{n}(-1)^{i}[g_0,\dots,\widehat{g_i},\dots,g_{n}]=0;
             $$
\item for any $A\in GL_n(K)$, we have $A\cdot[G]=[AG]$.
\end{enumerate}
Moreover, the relations (1) and (2) are consequences of relations (3) and (4).
\end{proposition}
\begin{proof}:
The proof of relations (1)-(5) are given in \cite[proposition 2.2]{AR}, so we will just prove the last statement. 
For simplicity, we write
$$
\partial[g_0,\cdots,g_n]=\sum_{i=0}^{n}(-1)^{i}[g_0,\dots,\widehat{g_i},\dots,g_{n}]=0.
$$
For the property (1), firstly, we notice that by (3) and (4):
$$
\partial[g_0,\cdots,g_{i+1},g_i, g_{i+1}\cdots,g_{n-1}]=
$$
$$(-1)^{i}[g_0,\dots,g_{i+1},g_i\dots,g_{n-1}]+(-1)^{i+2}[g_0,\dots,g_i,g_{i+1},\dots,g_{n-1}]=0$$
which means swapping of two vectors next to each other results a sign change. Therefore,
$$[g_0,\dots,g_i,\dots,g_j,\dots,g_{n-1}]=(-1)^{j-i}[g_0,\dots,g_j,g_i,\dots,g_{j-1},g_{j+1},\dots,g_{n-1}]$$
$$=(-1)^{j-i}(-1)^{j-i-1}[g_0,\dots,g_j,\dots,g_i,\dots,g_{n-1}].$$
For (2), it suffices to notice that by relations (3) and (4) we have
$$\partial[g_0,g_1,\cdots,g_{n-1}, ag_0]=[g_1,\dots,g_{n-1},ag_0]+(-1)^{n}[g_0,\dots,g_{n-1}]=0.$$
\end{proof}
\

\subsection{Some Homotopy Properties}
In this subsection, we generalize the results of \cite[Section 3]{KS}.
Let $X$ be a set and $E(X)$ be the simplicial complex defined at the beginning of section 4. Let $\q_f(X):=\Simpl(E(X))$ be the ordered set of non-empty finite subsets of $X$. By \cite[Lemma 1]{KS}, $\q_f(X)$ and $E(X)$ are contractible. 
\\

Let $V$ be a finite dimensional vector space over a field $K$ such that $\dim(V)\geq 1$. If $X\subset V$ is a non-empty finite subset, we denote $\langle X\rangle$ the subspace of $V$ generated by $X$.
\begin{definition}\label{def1}\
\begin{enumerate}
\item
We define $E^*(V)$ to be the subsimplicial complex of $E(V)$ such that $$Vert(E^*(V))=Vert(E(V))$$ and if $0\in X\in \Simpl(E^*(V))$ then $\langle X\rangle<V$.
\item
 $\q^*(V)$ is defined as the ordered set $\Simpl(E^*(V))$.
 \end{enumerate}
\end{definition}

\begin{remark}\label{rem4.5}
The simplicial complexes $E^*(V)$ and $E(V)$ are different from the ones defined in \cite{KS} since we need to include $0$ as a vertex here. In order to emphasize the difference between these two notations, we use $E^*(V-\{0\})$ and $E(V-\{0\})$ to denote the ones defined in loc. cit., i.e, $E(V-\{0\})$ is the simplicial complex of finite non-empty subsets of $V-\{0\}$ and $E^*(V-\{0\})\subset E(V-\{0\})$ is the sub simplicial complex of finite non-empty subsets $X$ such that $\langle X\rangle<V$.
\end{remark}
\

\begin{definition}\label{def2}
Let $V$ be a finite dimensional vector space and $W_0, W_1\subseteq V$ be two subspaces such that $W_0\leq W_1$. We define
$$
\q_f(W_0,W_1)=\{ X\in\q_f(W_1)\ |\ W_0\leq\langle X\rangle\}.
$$
\end{definition}

\begin{lemma}\label{lem0}
$
\q_f(W_0,W_1)
$
is contractible.
\end{lemma}
\begin{proof}
Let $B\subseteq W_1$ be a non-empty finite subset such that $W_0= \langle B\rangle$ and we define $\q_f(W_0,W_1)_B$ the sub ordered set of $\q_f(W_0,W_1)$ whose elements always contain $B$. Then the inclusion $\q_f(W_0,W_1)_B\hookrightarrow \q_f(W_0,W_1)$ has a left adjoint $X\mapsto X\cup\{B\}$. So $\q_f(W_0,W_1)_B$ and $\q_f(W_0,W_1)$ are homotopy equivalent. Since $\q_f(W_0,W_1)_B$ admits a smallest element $B$, it is contractible. It follows that $\q_f(W_0,W_1)$ is also contractible.

\end{proof}

Let $V$ be a vector space. Recall that (c.f, \cite[Section 2]{Q2}) $J(V)$ denotes the ordered set of pairs $(W_0,W_1)$ of $V$ such that $0\leq W_0\leq W_1\leq V$ and $\dim(W_1/W_0)<n$. The order is given by the relation that $(W_0,W_1)\leq (W'_0,W'_1)$ if and only if $W'_0\leq W_0\leq W_1\leq W'_1$. Such a pair $(W_0,W_1)\in J(V)$ is called a proper layer of $V$.
\begin{lemma}\label{lem1}
The order preserving map
$$
g': \Simpl(B\q^*(V))\to J(V);\ \ \ \ \ X_0<\cdots<X_p\mapsto(\langle X_0\rangle, \langle X_p\rangle)
$$
is a homotopy equivalence.
\end{lemma}
\begin{proof}
We will prove this lemma by using Quillen's theorem A.
Let $(W,W')\in J(V)$ be a proper layer of $V$. We consider the category
$$
g'\downarrow(W,W')=\{X_0<\cdots<X_p\in\Simpl(B\q^*(V))\ |\ W\leq\langle X_0\rangle\leq\langle X_p\rangle\leq W'\}.
$$
If $W'\neq V$, we see that $g'\downarrow(W,W')=\Simpl(B\q_f(W,W'))$ which is contractible by lemma \ref{lem0}.
\\

It leaves us to prove that $g'\downarrow(W,V)$ with $(W,V)\in J(V)$ is contractible. Let $\q^*(V)_W$ be the sub-ordered set of $\q^*(V)$ whose object $X$ satisfies $W\leq\langle X\rangle$. We denote $0\notin B\subset V$ a finite subset such that $W=\langle B\rangle$ and $\q^*(V)_B$ the sub-ordered set of $\q^*(V)$ whose objects contain $B$. The natural inclusion $\q^*(V)_B\hookrightarrow\q^*(V)_W$ admits a left adjoint $X\mapsto X\cup\{B\}$ and hence is a homotopy equivalence. Since $\q^*(V)_B$ has a minimal object $B$, it is contractible which implies that $\q^*(V)_W$ is also contractible. We conclude by noticing that $g'\downarrow(W,V)=\Simpl(B\q^*(V)_W)$ which is contractible. 

\end{proof}
\

\begin{proposition}[\cite{Q2}, section 2]\label{equiv1}
Suppose $n\geq 2$. There is a $GL(V)$-equivariant homotopy equivalence
$$
g: \Simpl(\Sigma T(V))\to J(V);\ \ \ \ \ (W_0<\cdots<W_p)\mapsto (W_0,W_p).
$$
\end{proposition}

\begin{corollary}\label{cor1'}
The simplicial complexes $E^*(V)$ and $\Sigma T(V)$ are homotopy equivalent.
\end{corollary}
\begin{proof}
According to proposition \ref{equiv1}, there is a homotopy equivalence 
$$
\Simpl(\Sigma T(V))\xrightarrow{\approx}J(V).
$$
Since $B\circ\Simpl=sd$, by lemma \ref{lem1}, we see that $\Simpl(sdE^*(V))$ and $\Simpl(\Sigma T(V))$ are homotopy equivalent. It follows that $sdE^*(V)$ and $\Sigma T(V)$ are homotopy equivalent (c.f, argument at the beginning of Chapter 4) and hence so do $E^*(V)$ and $\Sigma T(V)$.

\end{proof}

By Proposition \ref{equiv1} and Corollary \ref{cor1'}, we see that $E^*(V)$ and $J(V)$ are homotopy equivalent.
\\

Let $A$ be an integral domain with quotient field $K$ and $M$ be a torsion-free Noetherian  $A$-module. Write $V=K\otimes_A M$, so that $M$ is a lattice in $V$: we assume $\Dim(V)=n\ge 2$. A submodule $N$ of $M$ is \emph{pure} if $M/N$ is torsion-free. Let $G^*(M)$ be the poset of proper pure submodules of $M$ (those different from $0$ and $M$) and $G^*(V)$ be the poset of proper subspaces of $V$ (those different from $0$ and $V$). 
\begin{proposition}[\cite{K}, Prop. 4.2.4]\label{prop51}  
There exists mutually invertible bijections
$$-\otimes_AK: G^*(M)\longleftrightarrow G^*(V): -\cap M.$$
\end{proposition}
\

We define $Q^{\tf}(A)$ as Quillen's $Q$-construction on the category of finitely generated torsion-free $A$-modules and $Q_n^{\tf}(A)$ as the full subcategory of $Q^{\tf}(A)$ of torsion-free modules $M$ such that $rank(M)\leq n$.
\begin{lemma}\label{layer}
For any $d\in Q^{\tf}(A)$ of rank $n$, the functor $Q^{\tf}_{n-1}(A)\downarrow d\to J(d\otimes K)=J(V)$ sending
$$
\xymatrix{
 & c'\ar@{->>}[ld]_{\phi}\ar@{>->}[rd] & \\
 c & & d 
}
$$
to $(ker(\phi)\otimes K, c'\otimes K)$ is an isomorphism of categories.
\end{lemma}
\begin{proof}
Since $c'$ and $ker(\phi)$ are pure submodules of $d$, after proposition \ref{prop51}, we see that this functor is fully faithful and bijective on objects.  

\end{proof}
We notice that the inverse sends $(X,Y)\in J(V)$ to
$$
\xymatrix{
 & Y\cap d\ar@{->>}[ld]\ar@{>->}[rd] & \\
 Y\cap d/X\cap d & & d 
}
$$
\

\subsection{The Extended (Modular) Symbols}
Let $n:=\Dim(V)\geq 2$. By definition, we have $C_*(E(V))=C_*(E^*(V))$ for $*\leq n-1$. 
Since $E(V)$ is contractible (\cite[Lemma 1]{KS}), the complex $C_*(E(V))$ is acyclic for $*>0$. It follows that
$$
Ker(C_{n-1}(E^*(V))\xrightarrow{\partial}C_{n-2}(E^*(V)))
=Im(C_n(E(V))\xrightarrow{\partial}C_{n-1}(E(V)))
$$
So we have
$$
H_{n-1}(E^*(V))=Im(C_n(E(V))\xrightarrow{\partial}C_{n-1}(E(V)))/Im(C_n(E^*(V))\xrightarrow{\partial}C_{n-1}(E^*(V)))
$$
and the sequence
\begin{equation}\label{presentseq1}
C_{n+1}(E(V))\xrightarrow{\partial}C_n(E(V))/C_n(E^*(V))\xrightarrow{\partial} H_{n-1}(E^*(V))\to 0
\end{equation}
is exact.
By definition \ref{def1}, the symbol $(g_0,\cdots,g_i)\in C_i(E^*(V))$ if it satisfies one the following two conditions
\begin{itemize}
\item none of the vectors in $\{g_0,\cdots,g_i\}$ is zero.
\item the collection of vectors $\{g_0,\cdots,g_i\}$ contains zero and $\langle g_0,\cdots,g_i\rangle<V$.
\end{itemize}
Thus $C_n(E(V))/C_n(E^*(V))$ is a free Abelian group on symbols $(g_0,\cdots,g_n)$ such that one of the vectors is zero and the others are linearly independent. Moreover, if $g_0,\cdots,g_{n+1}\in V$, then 
$$
0=\partial\circ\partial(g_0,\cdots,g_{n+1})=\sum_{i=0}^{n+1}(-1)^i\partial(g_0,\cdots,\widehat{g_i},\cdots,g_{n+1})\in C_{n-1}(E(V))=C_{n-1}(E^*(V)).
$$ 
So, when $n=\Dim(V)\geq 2$, $H_{n-1}(E^*(V))$ is generated by symbols $\partial(g_0,g_1,\cdots,g_n)$ and presented by the relations
\begin{description}\label{present1}
\item[(a)] if the collection of vectors $\{g_0,\cdots,g_n\}$ does not contain zero or it contains zero but does not generate $V$, then $\partial(g_0,g_1,\cdots,g_n)=0$.
\item[(b)] the alternating sum 
          $$
          \sum_{i=0}^{n+1}(-1)^i\partial(g_0,\cdots,\widehat{g_i},\cdots,g_{n+1})
          $$ 
          equals zero.
\end{description}
\

We will call the symbols $\partial(g_0,\cdots,g_n)$ the extended (modular) symbols. 
\begin{proposition}\label{relationsms}\
The extended symbols $\partial(g_0,\cdots,g_n)$ satisfy the relations
\begin{description}
\item[(a0)] Swapping two vectors in $\partial(g_0, g_1, \cdots,g_n)$ changes the sign;
\item[(b0)] $\partial(ag_0, g_1,\cdots,g_n)=\partial(g_0,g_1,\cdots,g_n)$ for any $a\in K-\{0\}$.
\end{description}
\end{proposition}
\begin{proof}
We use the similar argument as given in the proof of proposition \ref{propms}.
\end{proof}
\

\begin{proposition}\label{presms}
When $n=\Dim(V)\geq 2$, The homology group $H_{n-1}(E^*(V))$ is generated by symbols $\partial(0,g_1,\cdots,g_n)$ and presented by the relations
\begin{enumerate}\label{present0}
\item If $g_1,\cdots, g_n$ are linearly dependent then $\partial(0,g_1,\cdots,g_n)=0$.
\item The alternating sum 
          $$
         \sum_{i=1}^{n+1}(-1)^i\partial(0,g_1,\cdots,\widehat{g_i},\cdots,g_{n+1})
          $$ 
          equals zero.
\end{enumerate}
\end{proposition}
\begin{proof}
We have seen that $H_{n-1}(E^*(V))$ is generated by symbols $\partial(g_0,g_1,\cdots,g_n)$ and presented by relations (a) and (b) above. In particular, these symbols satisfy the properties listed in proposition \ref{relationsms}. Since a non-zero symbol consists one zero vector, we can thus permute $0$ to the beginning after a proper sign change. It follows that $H_{n-1}(E^*(V))$ is generated by symbols $\partial(0,g_1,\cdots,g_n)$ which equal zero if $g_1,\cdots.g_n$ are linearly dependent, thus the relation 1.
\\

We notice that, by sequence (\ref{presentseq1}), if none of $g_0,\cdots,g_{n+1}$ is zero, then each summand in the alternating sum $\partial(g_0,\cdots,g_n)$ is zero in $C_n(E(V))/C_n(E^*(V))$. If two of these vectors are zero, the symbol $C_n(E(V))/C_n(E^*(V))\ni(g_0,\cdots, \widehat{g_i},\cdots,g_{n+1})$ is zero if $g_i\neq 0$. Suppose $g_j=g_k=0$, then by (a0) it is easy to see that 
$$
(-1)^j\partial(g_0,\cdots,\widehat{g_j},\cdots,g_n)+(-1)^k\partial(g_0,\cdots,\widehat{g_k},\cdots,g_n)=0.
$$ 
If more than three vectors are zero, we see that these vectors do not generate $V$, so by (a), each symbol $(g_0,\cdots, \widehat{g_i},\cdots,g_{n+1})=0$.
\\

It leaves us to consider the case where exactly one of these vectors is zero. 
According to the formula given in (b), if $g_0=0$, we get the relation 2. If, say, $g_j=0$ for some $0<j\leq n$, then we use (a0) of proposition \ref{relationsms} which implies that 
$$
\partial(g_0,\cdots,\widehat{g_i},\cdots,g_n)=\left\{
\begin{array}{ll}
(-1)^{j-1}\partial(0, g_0,\cdots,\widehat{g_i},\cdots, \widehat{g_j},\cdots,g_n); & i< j\\
(-1)^j\partial(0, g_0,\cdots, \widehat{g_j},\cdots, \widehat{g_i},\cdots,g_n); & i>j
\end{array}
\right.
$$
Moreover, by (a), we have $\partial(g_0,\cdots,\widehat{g_j},\cdots,g_n)=0$. It follows from (b) that
$$
\sum_{i=0}^{j-1}(-1)^{i+j-1}\partial(0, g_0,\cdots,\widehat{g_i},\cdots, \widehat{g_j},\cdots,g_n)+\sum_{i=j+1}^n(-1)^{i+j}\partial(0, g_0,\cdots, \widehat{g_j},\cdots, \widehat{g_i},\cdots,g_n)
$$
equals zero which implies the relation 2.

\end{proof}
\

\subsection{Modular Symbols and Extended Symbols}
For $n=\Dim(V)=1$, we have $H_0(E^*(V))=H_0(\Sigma T(V))=\Z\oplus\Z$. By corollary \ref{cor1'}, we have 
$$
\widetilde{St}(V)\simeq \widetilde{H}_{n-1}(E^*(V)),\ \ \ \ \ n\geq 1.
$$
In particular, if $n=1$, then $J(V)$ is discrete consisting of $(0,0)$ and $(V,V)$. We can apply our discussion in the previous section to this case and find that $\widetilde{H}_0(E^*(V))$ is generated by the symbols $\partial(0,g)$ with $\langle g\rangle=V$ which satisfy
\begin{enumerate}
\item $\partial(0,g)=\partial(g,0)$.
\item $\partial(0,g)=\partial(0,ag)$ for any non-zero scalar $a$.
\end{enumerate}
\

Suppose $\Dim(V)=n\geq 2$. Recall that (remark \ref{rem4.5}) $E^*(V-\{0\})$ be the simplicial complex of non-empty finite subsets $X$ of $V-\{0\}$ with $\Dim(\langle X\rangle)<n$. In \cite{KS}, we showed that there is a homotopy equivalence between $E^*(V-\{0\})$ and $T(V)$ (see \cite[theorem 1]{KS} and notice that the notations given there are different from those defined here in the sense of remark \ref{rem4.5}). 
The simplicial complexes $E^*(V)$ and $E^*(V-\{0\})$ are homotopy models (i.e, homotopy equivalent to) $\Sigma T(V)$ and $T(V)$ respectively (see corollary \ref{cor1'} and \cite[theorem 1]{KS}). 
\\

Consider the cocartesian diagram
$$
\xymatrix{
E^*(V-\{0\})\ar[r]^i\ar[d]_{i'} & E(V-\{0\})\ar[d]^j\\
\bigcup_{\Dim(W)=n-1}E(W)\ar[r]_{j'} & E^*(V)
}
$$
which gives a termwise short exact sequence
$$
0\to \widetilde{C}_*(E^*(V-\{0\}))\xrightarrow{(i_*,-i'_*)}\widetilde{C}_*(E(V-\{0\}))\oplus\widetilde{C}_*\left(\bigcup_WE(W)\right)
$$
$$
\xrightarrow{(j_*,j'_*)}\widetilde{C}_*(E^*(V))\to 0.
$$
Let $\underline{T}(V)$ be Quillen's model of the lower cone of $T(V)$ (c.f, \cite[p8]{Q2}) which is contractible.
We denote $P:=\bigcup_WE(W)$. Consider the functor
$$
AR: \Simpl(P)\to \underline{T}(V),\ \ \ \ \ \ Y\mapsto\langle Y\rangle.
$$
Since for every proper non-empty subspace (including $0$) $W'\subset V$, the comma category
$$
AR\downarrow W'=\Simpl(E(W'))
$$
is contractible, by Quillen's theorem A, $\Simpl(P)$ is homotopy equivalent to $\underline{T}(V)$. It follows that $\bigcup_WE(W)$ is contractible. 
Moreover, by \cite[Lemma 1]{KS}, $E(V-\{0\})$ is also contractible. 
Hence, the long exact sequence induced by the above short exact sequence gives an isomorphism
$$
\widetilde{H}_{n-1}(E^*(V))\to\widetilde{H}_{n-2}(E^*(V-\{0\})),\ \ \ \ \ \ \partial(0,g_1,\cdots,g_n)\mapsto\partial(g_1,\cdots,g_n).
$$ 
We define $\bar{C}_*$ to be the quotient complex fitting into the short exact sequence
\begin{equation}\label{seq5.3}
0\to C_*(E^*(V-\{0\}))\to C_*(E(V-\{0\}))\to \bar{C}_*\to 0
\end{equation}
which induces a resolution sequence (c.f, \cite[Section 6]{KS})
$$
\cdots\to\bar{C}_n\to\bar{C}_{n-1}\to H_{n-1}(\bar{C}_*)\to 0.
$$
Moreover, sequence (\ref{seq5.3}) gives an isomorphism
$$
H_{n-1}(\bar{C}_*)\xrightarrow{\simeq}\widetilde{H}_{n-2}(E^*(V-\{0\})),\ \ \ \ \ \overline{(g_1,\cdots,g_n)}\mapsto \partial(g_1,\cdots,g_n)
$$
where $\overline{(g_1,\cdots,g_n)}\in\bar{C}_{n-1}=C_{n-1}(E(V-\{0\}))/C_{n-1}(E^*(V-\{0\}))$ denotes a generator.
\\

\begin{theorem}\label{compare5}
When $n=\Dim(V)\geq 2$, the map
$$
\widetilde{H}_{n-1}(E^*(V))\longleftrightarrow \widetilde{St}(V);\ \ \ \ \ \partial(0,g_1,\cdots,g_n)\leftrightarrow[g_1,\cdots,g_n]
$$
is an isomorphism.
\end{theorem}
\begin{proof}
By our discussion, we have a series of isomorphisms
$$
\widetilde{H}_{n-1}(E^*(V))\simeq\widetilde{H}_{n-2}(E^*(V-\{0\}))\simeq H_{n-1}(\bar{C}_*)
$$
given by
$$
\partial(0,g_1,\cdots,g_n)\leftrightarrow\partial(g_1,\cdots,g_n)\leftrightarrow \overline{(g_1,\cdots,g_n)}.
$$
Moreover, \cite[Thm. 2]{KS} shows that there exists a diagram
$$
\xymatrix{
\cdots\ar[r] & \bar{C}_n\ar[r] & \bar{C}_{n-1}\ar[r]\ar[rd]_{ar} & H_{n-1}(\bar{C}_*)\ar[r]\ar@{-->}[d]^{\simeq} & 0\\
& & & \widetilde{St}(V)\ar[r] & 0
}
$$
such that $ar$ sends $\overline{(g_1,\cdots,g_n)}$ to the modular symbol $[g_1,\cdots,g_n]$. Hence the vertical map is an isomorphism sending $\overline{(g_1,\cdots,g_n)}$ to $[g_1,\cdots,g_n]$.

\end{proof}
\

Notice that there is also an isomorphism
$$
\widetilde{H}_{n-2}(E^*(V-\{0\}))\longleftrightarrow \widetilde{St}(V);\ \ \ \ \ \partial(g_1,\cdots,g_n)\leftrightarrow[g_1,\cdots,g_n].
$$
\

\section{Calculation of $d^1$ on Coefficents}
In this section, we will put $d^1$ into distinguished triangles of singular homology groups of simplicial complexes so that we can use the extended symbols (and hence modular symbols) to calculate $d^1$.
We focus on Quillen's case such that $A$ is an integral Noetherian domain and $\D=Q_n$ is the full sub-category of $Q^{\tf}(A)$ consisting of torsion-free modules of rank$\leq n$. Similarly, we define $\B=Q_{n-2}$ and $\C=Q_{n-1}$.
The functors $U:\B\hookrightarrow\C$ and $T:\C\hookrightarrow\D$ are natural inclusions which are connected cellular. 
According to lemma \ref{layer}, the category $\F_{T}(d)$ is isomorphic to the ordered set $J(V)$ of proper layers of $V$ so that we use a proper layer of $J(V)$ to denote an object in $\F_{T}(d)$ and vice versa. 
\\

\subsection{Some Homology Properties}
\begin{definition}[\cite{Sp}, p186]\label{MCdiag}
We say that a commutative diagram of simplicial complexes
$$
\xymatrix{
A\ar[r]^i\ar[d]_{i'} & B\ar[d]^j\\
C\ar[r]_{j'} & D
}
$$
is a Mayer-Vietoris diagram if all the maps are inclusions and 
$$D=B\cup C,\ \ \ \ \ A=B\cap C.
$$
\end{definition}
According to \cite[p186]{Sp}, a Mayer-Vietoris diagram induces a short exact sequence of (reduced) chain complexes
$$
0\to \widetilde{C}_*(A)\xrightarrow{(i_*,-i'_*)}\widetilde{C}_*(B)\oplus\widetilde{C}_*(C)\xrightarrow{(j_*,j'_*)}\widetilde{C}_*(D)\to 0
$$
where $\widetilde{C}_*$ denotes the reduced chain complex, i.e, $\widetilde{C_*}=Ker(C_*\to\Z)$.
\\

We define the simplicial complex
$$
E^{<n}(V):=\bigcup_{W<V, \Dim(W)=n-1}E(W).
$$  
Let $E^{n-1}(V)$ be the subset of $\Simpl (E^*(V))$ with elements $(0,X)$ such that $\Dim(\langle X\rangle)=n-1$.  We take $E^{**}(V):=E^*(V)-E^{n-1}(V)$. By definition \ref{defsc}, $E^{**}(V)$ is a sub-simplicial complexe of $E^*(V)$. 
Thus, we have a Mayer-Vietoris diagram
\begin{equation}\label{cocart2}
\xymatrix{
E^{<n}(V)\cap E^{**}(V)\ar[r]\ar[d] & E^{**}(V)\ar[d]\\
E^{<n}(V)\ar[r] & E^*(V)
}
\end{equation}
which induces a short exact sequence
$$
0\to\widetilde{C}_*(E^{<n}(V)\cap E^{**}(V))\to \widetilde{C}_*(E^{**}(V))\oplus \widetilde{C}_*(E^{<n}(V))
$$
\begin{equation}\label{exact1}
\to\widetilde{C}_*(E^*(V))\to 0.
\end{equation}
\

\begin{lemma}\label{cech}
Suppose that $(X,0)$ is a pointed simplicial complex ($0$ be the base point) which is covered by a family of pointed sub-simplicial complexes $\{X_i\}_{i\in I}$. If the intersection $\bigcap_{j\in J}X_j$ is contractible for all non-empty subsets $J\subseteq I$, then the natural map
$$
\bigvee_{i\in I}X_i\to X
$$
is a homology equivalence.
\end{lemma}
\begin{proof}
Let us define 
$$
X':=\bigvee_{i\in I}X_i
$$
which is also covered by $X'_i=X_i,\ i\in I$ such that $X'_i\cap X'_{i'}=0$ (the base point) for all $i\neq i'\in I$. So we have a map between $E^1$-terms of Mayer-Vietoris spectral sequences
$$
E^{'1}_{p,q}=H_p\left(\coprod_{|J|=q, j_i\neq j_k}X'_{j_0}\cap\cdots\cap X'_{j_q}\right)=\bigoplus_{|J|=q, j_i\neq j_k}H_p\left(X'_{j_0}\cap\cdots\cap X'_{j_q}\right)
$$
$$
\longrightarrow E^1_{p,q}=H_p\left(\coprod_{|J|=q, j_i\neq j_k}X_{j_0}\cap\cdots\cap X_{j_q}\right)=\bigoplus_{|J|=q, j_i\neq j_k}H_p\left(X_{j_0}\cap\cdots\cap X_{j_q}\right).
$$
Here, $|J|$ denotes the cardinality of the set $J$.
If the indices in $J\subseteq I$ are not identical then intersections of both sides are contractible and if the indices in $J$ are identical then the corresponding map $H_p(\bigcap_J X'_{j_k})\to H_p(\bigcap_JX_{j_k})$ is the identity map. So above map is an isomorphism. Since $E^1$ converges to $H_{p+q}\left(\bigcup_{i\in K}X_i\right)$ and $E^{'1}$ converges to $H_{p+q}\left(X'\right)$, we get an homology isomorphism
$$
H_{n}\left(X'\right)\xrightarrow{\simeq}H_{n}\left(\bigcup_{i\in K}X_i\right).
$$

\end{proof}
\

\begin{corollary}\label{equiv0}
$E^{<n}(V)\cap E^{**}(V)$ (resp. $E^{<n}(V)$) is homology equivalent to $\bigvee E^*(W)$ (resp. $\bigvee E(W)$) where the wedge sum is indexed by the set 
$$
\{W\ |\ W<V, \Dim(W)=n-1\}.
$$
\end{corollary}
\begin{proof}
In the category of simplicial complexes, we have
$$
E^{<n}(V)\cap E^{**}(V)=\bigcup_{W<V, \Dim(W)=n-1}E^*(W),\ \ \ \ \ \ \ \ \ E^{<n}(V)=\bigcup_{W<V, \Dim(W)=n-1}E(W)
$$
which are covered by $E^*(W)$ and $E(W)$ respectively and have $0$ as base points.
Since the $W$'s do not contain each other, any proper intersection
satisfies
$$
\bigcap_i E^*(W_i)=E\left(\bigcap_i W_i\right)
$$
and is contractible (\cite[Lemma 1]{KS}). It suffices to apply lemma \ref{cech}.

\end{proof}
This corollary implies that the distinguished triangle (\ref{exact1}) can be replaced by (up to quasi-isomorphism)
\begin{equation}\label{exact1'}
\bigoplus_{\Dim(W)=n-1}\widetilde{C}_*(E^*(W))\to \widetilde{C}_*(E^{**}(V))\to\widetilde{C}_*(E^*(V))\xrightarrow{+1}
\end{equation}
since $E^{<n}(V)$ has trivial reduced homology.
\\

\subsection[$d^1$ in Terms of Singular Homologies]{$d^1$ in Terms of Singular Homologies of Simplicial Complexes}
We compute (\ref{coeff}) in this subsection.
Since $\F_{T'}(d)$ is discrete and consists of admissible monomorphisms and admissible epimorphisms, we shall distinguish these two cases and get the formula of $d^1$ on coefficients by combining them.
\\

Recall that by the definition of relative reduced chains, we have
$$
\widetilde{C}_*(\F_T(d))=Ker(C_*(\F_T(d))\to \Z)=C_*(\widetilde{\F_T}(d)).
$$
So we will not distinguish these two notations in our discussion .
\\

\subsubsection[First Part of $d^1$]{The Formula for $d^1$ Having Image in Direct Sums of Reduced Steinberg Modules Indexed by Admissible Monomorphisms}
\

\paragraph*{Some Calculation of Singular Chains}\

\begin{lemma}
The simplicial map
$$
AR': sdE^*(V)\to \Sigma T(V),\ \ \ \ \ \ X\mapsto\langle X\rangle.
$$
is a homotopy equivalence.
\end{lemma}
\begin{proof}
According to lemma \ref{lem1}, the functor
$$
g': \Simpl(sdE^*(V))=\Simpl(B\q^*(V))\to J(V),\ \ \ \ \ (X_0<\cdots<X_p)\mapsto(\langle X_0\rangle, \langle X_p\rangle)
$$
is a homotopy equivalence.  So, combined with proposition \ref{equiv1}, we get the following commutative diagram
$$
\xymatrix{
\Simpl(sdE^*(V))\ar[rr]^{\Simpl\circ AR'}\ar[rrd]_{g'} & & \Simpl(\Sigma T(V))\ar[d]^g_{\approx}\\
 & & J(V)
}
$$
It follows that $\Simpl\circ AR'$ (and hence $B\circ\Simpl\circ AR'=sd\circ AR'$) is a homotopy equivalence. According to (\ref{e3}), for any simplicial complex $C$ there is a homotopy equivalence $\epsilon_C: |sdC|\xrightarrow{\approx}|C|$ which is natural on $C$. So, $|AR'|$ is a homotopy equivalence and so is $AR'$.  

\end{proof}
\

Moreover, for any sub simplicial complex $C\subseteq E^*(V)$, the map $AR'$ sends $sdC$ to a sub simplicial complex of $\Sigma T(V)$. In particular, $AR'(sdE^*(W))=\Sigma T(W)$ for $W$ a subspace of $V$.
Applying the functor $AR'\circ sd$ to (\ref{cocart2})  we get a Mayer-Vietoris diagram
\begin{equation}\label{cocart2'}
\xymatrix{
\bigcup_{\Dim(W)=n-1}\Sigma T(W)\ar[rr]\ar[d] & & AR'(sdE^{**}(V))\ar[d]\\
AR'(sdE^{<n}(V))\ar[rr] & & \Sigma T(V)
}
\end{equation}
where $AR'(sdE^{<n}(V))=\bigcup_WCt(W)$ such that $Ct(W)$ is the simplicial complex (c.f, \cite[Section 2]{Q2}) whose $p$-simplices are
$$
W_0<\cdots< W_p
$$
with $W_i$ subspaces (may equals $0$ or $W$) for each $0\leq i\leq p$. It follows that, by lemma \ref{cech}, $AR'(sdE^{<n}(V))$ has trivial reduced homology since the $Ct(W)$'s and their intersections are contractible.
Lemma \ref{cech} also gives us the following canonical quasi-isomorphisms
$$
\bigoplus_{\Dim(W)=n-1}\widetilde{C}_*(E^*(W))\xrightarrow{\sim} \widetilde{C}_*\left(\bigcup_{\Dim(W)=n-1}E^*(W)\right)
$$
and
$$
\bigoplus_{\Dim(W)=n-1}\widetilde{C}_*(\Sigma T(W))\xrightarrow{\sim} \widetilde{C}_*\left(\bigcup_{\Dim(W)=n-1}\Sigma T(W)\right).
$$
Thus,
the functor $AR'\circ sd$ from diagram (\ref{cocart2}) to diagram (\ref{cocart2'}) gives a morphism between distinguished triangles
\begin{equation}\label{cocart2''}
\xymatrix{
\bigoplus_{\Dim(W)=n-1}\widetilde{C}_*(E^*(W))\ar[r]\ar[d] & \widetilde{C}_*(E^{**}(V))\ar[r]\ar[d] & \widetilde{C}_*(E^*(V))\ar[r]^{\ \ \ \ \ \ \ +1}\ar[d] & \\
\bigoplus_{\Dim(W)=n-1}\widetilde{C}_*(\Sigma T(W))\ar[r] & \widetilde{C}_*(AR'(sdE^{**}(V)))\ar[r] & \widetilde{C}_*(\Sigma T(V))\ar[r]^{\ \ \ \ \ \ \ +1} &
}
\end{equation}
such that all vertical morphisms are quasi-isomorphisms.
\\

\paragraph*{$d^1$ in Terms of the Boundary Map of Mayer-Vietoris Sequence}\

Recall that there is a homotopy equivalence (c.f, \cite[Prop. p10]{Q2})
$$
g: \Simpl(\Sigma T(V))\to J(V),\ \ \ \ \ \ (W_0<\cdots<W_p)\mapsto(W_0,W_p).
$$
We define $J_1:=g\circ\Simpl(AR'(sdE^{**}(V)))\subset J(V)$ which is an ordered set consisting of proper layers $(X,Y)$ such that if $Y\neq V$ then $\Dim(Y/X)\leq n-2$ and $(L,V)$ with $\Dim(L)=1$. In particular, the above order preserving map $g$ gives a homotopy equivalence between $\Simpl(AR'(sdE^{**}(V)))$ and $J_1$. Therefore, we obtain a morphism of distinguished triangles
\begin{equation}\label{cocart2'''}
\xymatrix{
\bigoplus_{\Dim(W)=n-1}\widetilde{C}_*(\Sigma T(W))\ar[r]\ar[d] & \widetilde{C}_*(AR'(sdE^{**}(V)))\ar[r] \ar[d]& \widetilde{C}_*(\Sigma T(V))\ar[r]^{\ \ \ \ \ \ \ +1}\ar[d] & \\
\bigoplus_{\Dim(W)=n-1}\widetilde{C}_*(J(W))\ar[r] & \widetilde{C}_*(J_1)\ar[r] & \widetilde{C}_*(J(V))\ar[r]_{\ \ \ \ \ \ +1} &
}
\end{equation}
such that all vertical morphisms are induced by the functor $g\circ\Simpl$ and hence are quasi-isomorphisms.
\\

For any $d\in\D-\C$ such that $V=d\otimes K$, the inclusion functor $i': \F_{TU}(d)\hookrightarrow J_1$ is connected cellular which gives rise to a homotopy cocartesian diagram
$$
\xymatrix{
(J_1-\F_{TU}(d))\int\F_{i'}\ar[rr]\ar[d] & & \F_{TU}(d)\ar[d]^{i'}\\
J_1-\F_{TU}(d)\ar[rr] & & J_1
}
$$
Since $\F_{TU}(d)\subset J(V)$ consists of proper layers $(X,Y)$ such that $\Dim(Y/X)\leq n-2$, by our description of $J_1$, we see that $J_1-\F_{TU}(d)$ is a set (i.e, discrete category) consisting of proper layers $(L,V)$ with $\Dim(L)=1$. Since $(L,V)\int\F_{i'}\simeq\F_{U_*}(L,V)$, 
the homotopy cocartesian diagram above gives us a distinguished triangle
\begin{equation}\label{dist1}
\bigoplus_{d\twoheadrightarrow c}C_*(\widetilde{\F_U}(c))\to C_*(\widetilde{\F_{TU}}(d))\to \widetilde{C}_*(J_1)\xrightarrow{+1}
\end{equation}
Moreover, the lower triangle of (\ref{cocart2'''}) fits into the commutative diagram of distinguished triangles
\begin{equation}\label{part1}
\xymatrix{
\bigoplus_{d\twoheadrightarrow c\in\F_{T'}(d)}C_*(\widetilde{\F_U}(c))\ar[d]\ar[r]^= & \bigoplus_{d\twoheadrightarrow c}C_*(\widetilde{\F_U}(c))\ar[d] & & \\
\bigoplus_{c\to d\in\F_{T'}(d)}C_*(\widetilde{\F_U}(c))\ar[r]\ar[d] & C_*(\widetilde{\F_{TU}}(d))\ar[r]\ar[d] & C_*(\widetilde{\F_T}(d))\ar[d]^{=}\ar[r]^{\ \ \ \ \ \ d^1} &  \\
\bigoplus_{\Dim(W)=n-1}\widetilde{C}_*(J(W))\ar[r] & \widetilde{C}_*(J_1)\ar[r] & \widetilde{C}_*(J(V))\ar[r]_{\ \ \ \ \ \ \ +1} & 
}
\end{equation}
where the left vertical sequence is split short exact, the middle vertical distinguished triangle is the triangle (\ref{dist1}) and the middle horizontal one is the triangle (\ref{IIc'''''}). 
Hence, there exists a commutative diagram
$$
\xymatrix{
H_{n-1}(\widetilde{\F_T}(d))\ar[rrr]^{d^1}\ar[d]^= & & & \bigoplus_{c\to d\in\F_{T'}(d)}H_{n-2}(\widetilde{\F_U}(c))\ar[d]\\
\widetilde{H}_{n-1}(J(V))\ar[rrr]_{+1} & & & \bigoplus_{\Dim(W)=n-1}\widetilde{H}_{n-2}(J(W))
}
$$
The lower map of the above diagram is calculated as the boundary map of the lower triangle of diagram (\ref{part1}) which gives part of $d^1$. We will give its formula in this subsection.
\\

\paragraph*{The Formula}\

Together with (\ref{cocart2''}) and (\ref{cocart2'''}), up to explicit quasi-isomorphisms, we are able to calculate one part of $d^1$ on coefficients through the boundary map of (\ref{exact1'}) and hence of (\ref{exact1}) which fits into a Mayer-Vietoris sequence after applying the functor $H_*(-)$. More precisely, we take $[a]\in\widetilde{Z}_{n-1}(E^*(V))$ and lift it to $[a]=[a_1]+[a_2]\in \widetilde{C}_{n-1}(E^{<n}(V))\oplus\widetilde{C}_{n-1}(E^{**}(V))$. Then we apply the differential map to $[a_2]$ to get $\partial([a_2])=-\partial([a_1])$. The boundary map is then given by $[a]\mapsto-\partial([a_1])$.
\\

We write $\widetilde{H}_{n-1}(E^*(V))=\widetilde{St}(d)$ and $\widetilde{H}_{n-2}(E^*(W))=\widetilde{St}(c)$ for $d\in\D-\C$ of rank $n$ with $V=d\otimes K$ and $c\subset d$ a sub-module of rank $n-1$ with $W=c\otimes K$.
\begin{proposition}\label{prop2}
The boundary map of (\ref{exact1}) is given by
$$
\widetilde{St}(d)\to\bigoplus_{c\rightarrowtail d\in\F_{T'}(d)}\widetilde{St}(c),\ \ \ \ \ \ \partial(0,g_1,\cdots,g_n)\mapsto -\sum_{i=1}^n(-1)^i((0,W_i), \partial(0,g_1,\cdots,\widehat{g_i},\cdots, g_n))
$$
where $d\in\D-\C$, $W_i:=\langle g_1,\cdots,\widehat{g_i},\cdots, g_n\rangle$ for $1\leq i\leq n$ and $(0,W_i)$ stands for the index of the summand.
\end{proposition}
\begin{proof}
We will use extended symbols. Consider
$$
[a]:=\partial(0,g_1,\cdots,g_n)=\sum_{i=1}^n(-1)^i(0,g_1,\cdots,\widehat{g_i},\cdots,g_n)+(g_1,\cdots,g_n)\in \widetilde{Z}_{n-1}(E^*(V))
$$
such that
$$
(0, g_1,\cdots,\widehat{g_i},\cdots,g_n)\in C_{n-1}(E(W_i)),\ \ \ \ \  (g_1,\cdots,g_n)\in C_{n-1}(E^{**}(V)).
$$
It follows that we can lift $[a]$ to 
$$
[a']=[a_1]+[a_2]=\sum_{i=1}^n(-1)^i(0,g_1,\cdots,\widehat{g_i},\cdots,g_n)+(g_1,\cdots,g_n)\in C_*(E^{<n}(V))\oplus C_*(E^{**}(V))
$$
where $[a_1]=\sum(-1)^i(0,g_1,\cdots,\widehat{g_i},\cdots,g_n)\in C_{n-1}(E^{<n}(V))$ and $[a_2]=(g_1,\cdots,g_n)\in C_{n-1}(E^{**}(V))$.
Apply $\partial$ to $[a_2]$, we have
$$
\partial[a_2]=-\partial[a_1]=-\sum_{i=1}^n(-1)^i\partial(0,g_1,\cdots,\widehat{g_i},\cdots,g_n).
$$
So we get the formula we are looking for
$$
\partial(0,g_1,\cdots,g_n)\longmapsto -\sum_{i=1}^n(-1)^i((0,W_i), \partial(0,g_1,\cdots,\widehat{g_i},\cdots,g_n)).
$$

\end{proof}
\

\subsubsection[Second Part of $d^1$]{The Formula for $d^1$ Having Image in Direct Sums of Reduced Steinberg Modules Indexed by Admissible Epimorphisms}
In this section, we will calculate the other part of $d^1$ than proposition \ref{prop2}.
\\

\paragraph*{Some distinguished Triangles}\

Suppose $\Dim(V)\geq 2$. We define $J_2$ as the sub-ordered set of $J(V)$ obtained by removing proper layers $(L,V)$ with $\Dim(L)=1$. Thus, $J(V)-J_2=\{(L,V)\ |\ \Dim(L)=1\}$ and the inclusion $i'':J_2\hookrightarrow J(V)$ is connected cellular. We have the following homotopy cocartesian diagram
$$
\xymatrix{
\coprod_{(L,V)}J(L,V)\ar[r]\ar[d] & J_2\ar[d]^{i''}\\
\coprod_{(L,V)}\star\ar[r] & J(V)
}
$$
where $J(L,V)=i''\downarrow(L,V)$ denotes the ordered set of proper layers $\{(W,W')\ | (W,W')<(L,V)\}$. This diagram induces a distinguished triangle
$$
\bigoplus_{\Dim(L)=1}\widetilde{C}_*(J(L,V))\to \widetilde{C}_*(J_2)\to \widetilde{C}_*(J(V))\xrightarrow{+1}
$$
\

We write $X$ as the sub-simplicial complex of $B\q^*(V)=sdE^*(V)$ obtained by deleting the simplices 
$$
\{K_0<\cdots<K_p\ |\ \Dim(\langle K_0\rangle)=1,\ \langle K_p\rangle=V\}.
$$
\begin{lemma}
The order preserving map
$$
g': \Simpl(X)\to J_2,\ \ \ \ \ \  (K_0<\cdots<K_p)\mapsto (\langle K_0\rangle, \langle K_p\rangle)
$$
is a homotopy equivalence.
\end{lemma}
\begin{proof}
It suffices to restrict $g'$ of lemma \ref{lem1} to $\Simpl(X)$.

\end{proof}
\

Suppose $\Dim(L)=1$. Let $B\q(V-\{0\})_L$ be the sub-simplicial complex of $B\q(V-\{0\})=sdE(V-\{0\})$ whose vertices $K$ satisfy $L\leq\langle K\rangle$ and we define $X_L$ as the sub-simplicial complex of $B\q(V-\{0\})_L$ obtained by deleting the set of simplices
$$
\{K_0<\cdots<K_p\ |\ \langle K_0\rangle=L,\ \langle K_p\rangle=V\}.
$$ 
For two subspaces $W\subseteq W'$ of $V$, we denote 
$$
\q'_f(W,W')=\{K\in\q_f(W'-\{0\})\ |\ W\leq\langle K\rangle\}.
$$
\begin{lemma}\label{l49}\

\begin{enumerate}
\item The ordered set $\q'_f(W,W')$ is contractible.
\item The functor $g'$ restricted on $X_L$ induces a homotopy equivalence 
$$
g': \Simpl(X_L)\to J(L,V).
$$
\end{enumerate}
\end{lemma}
\begin{proof}\

\begin{enumerate}

\item
Let $A\subset W'-\{0\}$ be a finite subset such that $\langle A\rangle=W$ and $\q'_f(W,W')_A$ be the sub-simplicial complex of $\q'_f(W,W')$ whose elements contain $A$. The natural inclusion $\q'_f(W,W')_A\hookrightarrow\q_f(W,W')$ admits a left adjoint $B\mapsto B\cup A$, so the two ordered sets are homotopy equivalent. Since $\q'_f(W,W')_A$ has a minimal element $A$, it is contractible. It follows that $\q'_f(W,W')$ is contractible. 

\item
For any proper layer $(W,W')<(L,V)$, we consider the comma category
$$
g'\downarrow(W,W')=
$$
$$\{K_0\leq\cdots\leq K_p\in\Simpl(B\q_f(V-\{0\}))\ |\ W\leq\langle K_0\rangle\leq\langle K_p\rangle\leq W'-\{0\}\}
$$
$$
=\Simpl(B\q'_f(W, W')).
$$
We conclude by $1.$ and Quillen's theorem A.
\end{enumerate}

\end{proof}
\

\begin{corollary}\label{c50}
The natural chain map 
$$
\bigoplus_{\Dim(L)=1}\widetilde{C}_*\left(X_L\right)\to \widetilde{C}_*\left(\bigcup_{\Dim(L)=1}X_L\right)
$$
is a quasi-isomorphism.
\end{corollary}
\begin{proof}
By the proof of the previous lemma, it is easy to show that $\bigcup_L sdX_L\to \bigcup_L BJ(L,V)$ is a homotopy equivalence by using Quillen's theorem A. Moreover, $BJ(L,V)$ is a simplicial complex pointed at $(V,V)$. If $L\neq L'\subset V$ we have $BJ(L,V)\cap BJ(L',V)=BJ'(L\oplus L',V)$ where $J'(L\oplus L',V)$ denotes the ordered set consisting of proper layers $(W,W')\leq (L\oplus L',V)$. Since $J'(L\oplus L',V)$ has a maximal element $(L\oplus L',V)$ it is contractible. Similarly, we can prove that $\bigcap_{L\in S}BJ(L,V)$ is also contractible for any subset $S$ of $\{L\ |\ \Dim(L)=1\}$.
Then, lemma \ref{cech} shows that $\bigvee_LBJ(L,V)\to \bigcup_LBJ(L,V)$ is a homology equivalence. Thus, we have a commutative diagram
$$
\xymatrix{
\bigoplus_L\widetilde{C}_*(sdX_L)\ar@{-->}[r]\ar[d] & \widetilde{C}_*\left(\bigcup_L sdX_L\right)\ar[d]\\
\bigoplus_L\widetilde{C}_*(BJ(L,V))\ar[r] & \widetilde{C}_*\left(\bigcup_LBJ(L,V)\right)
}
$$
where all solid morphisms are quasi-isomorphisms, it follows that the dotted one is so. Hence the map in question is a quasi-isomorphism.

\end{proof}
\

Since the diagram
$$
\xymatrix{
\bigcup_{\Dim(L)=1}X_L\ar[r]\ar[d] & X\ar[d]\\
sdE(V-\{0\})\ar[r] & sdE^*(V)
}
$$
is Mayer-Vietoris, together with corollary \ref{c50},
we get a short exact sequence fitting into the commutative diagram of distinguished triangles
$$
\xymatrix{
\bigoplus_{\Dim(L)=1}\widetilde{C}_*\left(X_L\right)\ar[r]\ar[d] & \widetilde{C}_*(X)\oplus\widetilde{C}_*(sdE(V-\{0\}))\ar[r]\ar[d]^{(g'\circ\Simpl,0)} & \widetilde{C}_*(sdE^*(V))\ar[d]\\
\bigoplus_{\Dim(L)=1}\widetilde{C}_*(J(L,V))\ar[r] &  \widetilde{C}_*(J_2)\ar[r] &  \widetilde{C}_*(J(V)) 
}
$$
\begin{equation}\label{exact3}
\xymatrix{
\ar[r] & \bigoplus_{L}\widetilde{C}_*\left(X_L\right)[1]\ar[d] \\
\ar[r] & \bigoplus_{L}\widetilde{C}_*(J(L,V))[1]
}
\end{equation}
where all vertical morphisms are quasi-isomorphisms.
\\

As mentioned at the beginning of this section, we identify $\F_{T'}(d)$ with $J(V)$ and use a proper layer to denote an object of $\F_{T'}(d)$ and vice versa.  Notice that $J_2$ consists of proper layers $(W,W')$ with $\Dim(W'/W)\leq n-2$ and $(0,W)$ with $\Dim(W)=n-1$. Thus $\F_{TU}(d)\hookrightarrow J_2$ is connected cellular and we have a homotopy cocartesian diagram
$$
\xymatrix{
\coprod_{(0,W)}\F_{U}(W)\simeq\coprod_{(0,W)}\F_{U_*}(0,W)\ar[rr]\ar[d] & & \F_{TU}(d)\ar@{^(->}[d]\\
\coprod_{(0,W)}\star\ar[rr] & & J_2
}
$$
which gives rise to a distinguished triangle
\begin{equation}\label{dist2}
\bigoplus_{c\rightarrowtail d\in\F_{T'}(d)}C_*(\widetilde{\F_U}(c))\to C_*(\widetilde{\F_{TU}}(d))\to  C_*(J_2)\xrightarrow{+1}
\end{equation}
Moreover, the lower triangle of (\ref{exact3}) fits into the commutative diagram of distinguished triangles
\begin{equation}\label{part2}
\xymatrix{
\bigoplus_{c\rightarrowtail d\in\F_{T'}(d)}C_*(\widetilde{\F_U}(c))\ar[d]\ar[r]^= & \bigoplus_{c\rightarrowtail d\in\F_{T'}(d)}C_*(\widetilde{\F_U}(c))\ar[d] & & \\
\bigoplus_{c\to d\in\F_{T'}(d)}C_*(\widetilde{\F_U}(c))\ar[r]\ar[d] & C_*(\widetilde{\F_{TU}}(d))\ar[r]\ar[d] & C_*(\widetilde{\F_T}(d))\ar[d]^{=}\ar[r]^{\ \ \ \ \ \ d^1} &  \\
\bigoplus_{\Dim(L)=1}\widetilde{C}_*(J(L,V))\ar[r] & \widetilde{C}_*(J_2)\ar[r] & \widetilde{C}_*(J(V))\ar[r]_{\ \ \ \ \ \ \ +1} & 
}
\end{equation}
where the left vertical sequence is split short exact, the middle vertical distinguished triangle is the triangle (\ref{dist2}) and the middle horizontal one is the triangle (\ref{IIc'''''}). 
\\

\paragraph*{$d^1$ in Terms of the Boundary Map of Mayer-Vietoris Sequence}\

Our discussion above implies that we have a commutative diagram
$$
\xymatrix{
H_{n-1}(\widetilde{\F_T}(d))\ar[rrr]^{d^1}\ar[d]^= & & & \bigoplus_{c\to d\in\F_{T'}(d)}H_{n-2}(\widetilde{\F_U}(c))\ar[d]\\ai
\widetilde{H}_{n-1}(J(V))\ar[rrr]_{+1} & & & \bigoplus_{\Dim(L)=1}\widetilde{H}_{n-2}(J(L,V))
}
$$
The lower map of the above diagram is calculated as the boundary map of the lower triangle of diagram (\ref{part2}) which gives the other part of $d^1$. We will give its formula in this subsection.
So, together with (\ref{exact3}), we see that the remaining part of $d^1$ can be calculated via the boundary map of the top distinguished triangle of (\ref{exact3}). 
\\

Let $X'_L:=sdE^*(V/L)$.  The quotient map $V-\{0\}\to V/L$ induces a simplicial map $q_L: X_L\to X'_L$.  By lemma \ref{lem1}, the functor 
$$
g': \Simpl(X'_L)\to J(V/L),\ \ \ \ \ \ (K'_0,\cdots,K'_p)\mapsto(\langle K'_0\rangle, \langle K'_p\rangle)
$$
is a homotopy equivalence. Moreover, the functor
$$
J(L,V)\to J(V/L),\ \ \ \ \ (W,W')\mapsto (W/L=\overline{W},W'/L=\overline{W'})
$$
is an isomorphism.
So there exists a commutative diagram
$$
\xymatrix{
X_L\ar[rr]^{q_L}\ar[d]_{\approx}^{g'} & & X'_L\ar[d]^{\approx}_{g'}\\
J(L,V)\ar[rr]_{\simeq} & & J(V/L)
}
$$
It follows that $q_L$ is a homotopy equivalence. 
Moreover, in the following diagram 
$$
\xymatrix{
\bigoplus_{\Dim(L)=1}\widetilde{C}_*\left(X_L\right)\ar[r]\ar[d]^{\bigoplus q_{L*}} & \widetilde{C}_*(X)\oplus\widetilde{C}_*(sdE(V-\{0\}))\ar[r]\ar[d]^{=} & \widetilde{C}_*(sdE^*(V))\ar[d]^=  \\
\bigoplus_{\Dim(L)=1}\widetilde{C}_*\left(sdE^*(V/L)\right)\ar[r]\ar[d]^{\epsilon_*} &  \widetilde{C}_*(X)\oplus\widetilde{C}_*(sdE(V-\{0\}))\ar[r]\ar[d]^{(1,\epsilon_*)} & \widetilde{C}_*(sdE^*(V))\ar[d]^{\epsilon_*} \\
\bigoplus_{\Dim(L)=1}\widetilde{C}_*\left(E^*(V/L)\right)\ar[r] &  \widetilde{C}_*(X)\oplus\widetilde{C}_*(E(V-\{0\}))\ar[r] & \widetilde{C}_*(E^*(V))
}
$$
\begin{equation}\label{exact3'}
\xymatrix{
\ar[r] & \bigoplus_{\Dim(L)=1}\widetilde{C}_*\left(X_L\right)[1]\ar[d] \\
\ar[r] & \bigoplus_{\Dim(L)=1}\widetilde{C}_*\left(sdE^*(V/L)\right)[1]\ar[d]\\
\ar[r] & \bigoplus_{\Dim(L)=1}\widetilde{C}_*\left(E^*(V/L)\right)[1]
}
\end{equation}
the middle and lower triangles are quasi-isomorphic to the top one, so they are distinguished triangles.
It suggests that up to explicit quasi-isomorphisms we can use the bottom triangle to calculate our formula. By the naturality of $\epsilon$ of (\ref{e3}), we first calculate the boundary map of the second horizontal distinguished triangle. Since the top triangle is Mayer-Vietoris, we take $[z]\in\widetilde{Z}_{n-1}(sdE^*(V))$ and lift it as $[z]=[z_1]+[z_2]$ such that $[z_1]\in\widetilde{C}_{n-1}(X)$ and $[z_2]\in\widetilde{C}_{n-1}(sdE(V-\{0\}))$. Then we apply the differential map to get the homology class
$$
\partial([z_1])=-\partial([z_2])\in\widetilde{H}_{n-2}\left(\bigcup_LX_L\right)\simeq\bigoplus_L\widetilde{H}_{n-2}(X_L).
$$ 
We need to calculate the image of $\partial([z_1])=-\partial([z_2])$ in $\bigoplus\widetilde{H}_{n-2}(sdE^*(V/L))$ under $\bigoplus q_{L*}$.
\\

There exists a commutative diagram for each $L\subset V$ with $\dim(L)=1$
$$
\xymatrix{
\widetilde{C}_{n-1}\left(\bigcup_L B\q(V-\{0\})_L\right)\ar@{=}[r] & \widetilde{C}_{n-1}(sdE(V-\{0\}))\ar[d]^{\partial}\ar[r]^{q'_L} &
 \widetilde{C}_{n-1}(sdE(V/L))\ar[d]^{\partial}\\
\widetilde{Z}_{n-2}\left(\bigcup_LX_L\right)\ar@{^(->}[r] & \widetilde{Z}_{n-2}(sdE(V-\{0\}))\ar[r]_{q'_L} &
 \widetilde{Z}_{n-2}(sdE^*(V/L))
}
$$
where $q'_L$ is induced by the simplicial map $E(V-\{0\})\to E(V/L)$ sending a vertex $v$ to the quotient class $\bar{v}$ represented by $v$.
We write 
$$
[z_2]=\sum_L n_L[y_L],\ \ \ \ \ [y_L]\in \widetilde{C}_{n-1}(B\q(V-\{0\})_L).
$$
In particular, if $[y_L]$ also belong to $\widetilde{C}_{n-1}(B\q(V-\{0\})_{L'})$ with $L\neq L'$, then it must lie in $\widetilde{C}_{n-1}(X_L)$ (and in $\widetilde{C}_{n-1}(X_{L'})$) and hence
$$
0=\partial[y_L]\in\widetilde{H}_{n-2}\left(\bigcup_L X_L\right).
$$
Therefore, we may write
$$
\partial[z_1]=-\partial[z_2]=\sum n_L\partial[y_L]\in\bigoplus_L\widetilde{H}_{n-2}(X_L)
$$ 
such that each $[y_L]$ in a unique $\widetilde{C}_{n-1}(B\q(V-\{0\})_L)$.
Moreover, if $[y_{L'}]\in\widetilde{C}_{n-1}(B\q(V-\{0\})_{L'})$ and $L\neq L'\subset V$, then 
$$
0=q'_L(\partial[y_{L'}])=\partial q'_L([y_{L'}])\in\widetilde{H}_{n-2}(sdE^*(V/L))
$$ 
since $q'_L[y_{L'}]\in\widetilde{C}_{n-1}(sdE^*(V/L))$.
Thus we find that
$$
\sum_Lq'_L\left(\sum n_L\partial[y_L]\right)=\sum n_L q'_L(\partial[y_L])\in\bigoplus_L\widetilde{H}_{n-2}(sdE^*(V/L))
$$
and it suffices to calculate the image of $-\partial([z_2])$ under the morphisms
$$
\widetilde{Z}_{n-2}(sdE(V-\{0\}))\xrightarrow{q'_L}\widetilde{Z}_{n-2}(sdE^*(V/L))\to \widetilde{H}_{n-2}(sdE^*(V/L))
$$
for all lines $L\subset V$.
\\

\begin{proposition}\label{prop2'}
We write $\widetilde{St}(V/L):=\widetilde{H}_{n-2}(E^*(V/L))$. The second part of $d^1$ on coefficients 
$$
\widetilde{St}(d)\to\bigoplus_{(L,V)\in\F_{T'}(d)}\widetilde{St}(V/L)
$$
is given by
$$
\partial(0,g_1,\cdots,g_n)\mapsto \sum_{i=1}^n(-1)^i((L_i,V), \partial(0, \bar{g}_1, \cdots, \widehat{g_i},\cdots,\bar{g}_n))
$$
where $L_i:=\langle g_i\rangle$ and $\bar{g}_j$ denotes the image of $g_j$ under the quotient map $V\to V/L_i$ for all $1\leq i\neq j\leq n$. The proper layers $(L_i,V)$ stand for the index of the summands.  
\end{proposition}
\begin{proof}
By our discussion above and the naturality of $\epsilon$ of (\ref{e3}), the boundary map of the bottom triangle can be calculated as follows. Let $[a]\in\widetilde{Z}_{n-1}(E^*(V))$ and we lift $[a]$ to 
$$
[a']=[a'_1]+[a'_2]\in  \widetilde{C}_{n-1}(X)\oplus\widetilde{C}_{n-1}(E(V-\{0\})).
$$ 
We have $\partial([a'_1])\in\widetilde{Z}_{n-2}(X)$ and $\partial([a'_2])\in\widetilde{Z}_{n-2}(E(V-\{0\}))$. Then, we calculate the image of $-\partial([a'_2])$ under the morphism 
$$
\widetilde{Z}_{n-2}(E(V-\{0\}))\to \widetilde{Z}_{n-2}(E^*(V/L))\to \widetilde{H}_{n-2}(E^*(V/L))
$$
to find the image of $[a]$ in $\bigoplus_L\widetilde{St}(V/L)$. Here, the first map is induced by the simplicial map $E(V-\{0\})\to E(V/L)$ sending a vertex $v$ to the quotient class $\bar{v}$ represented by $v$.
\\

Let us take
$$
0\neq[a]=\partial(0,g_1,\cdots,g_n)=\sum_{i=1}^n(-1)^i(0,g_1,\cdots,\widehat{g_i},\cdots,g_n)+(g_1,\cdots,g_n)\in \widetilde{Z}_{n-1}(E^*(V))
$$
and write $[a]=[a_1]+[a_2]$ such that
$$
[a_1]=\sum_{i=1}^n(-1)^i(0,g_1,\cdots,\widehat{g_i},\cdots,g_n),\ \ \ \ \ \ \ \ \  [a_2]=(g_1,\cdots,g_n).
$$
We lift $[a]=[a_1]+[a_2]$ to $[a']=[a'_1]+[a'_2]$ with $[a'_1]\in\widetilde{C}_*(X)$ and 
$$
[a'_2]=[a_2]=(g_1,\cdots,g_n)\in\widetilde{C}_*(E(V-\{0\})).
$$
By our discussion in section 4.3, $\widetilde{St}(V/L)$ is generated by the extended symbols $\partial(\bar{v}_0,\cdots,\bar{v}_{n-1})$ such that if none of $\bar{v}_0,\cdots,\bar{v}_{n-1}$ is zero or some $\bar{v}_i=0$ but they do not span $V/L$ then 
$\partial(\bar{v}_0,\cdots,\bar{v}_{n-1})=0$. So the image of $-\partial([a_2])$ under the morphism
$$
\widetilde{Z}_{n-2}(E(V-\{0\}))\to\widetilde{H}_{n-2}(E^*(V/L))=\widetilde{St}(V/L)
$$
equals zero if none of $g_i, 1\leq i\leq n$ generates $L$, and equals 
$$
-\partial(\bar{g}_1,\cdots, \bar{g}_{i-1},0,\bar{g}_{i+1},\cdots, \bar{g}_n)=(-1)^i(0,\bar{g}_1,\cdots,\widehat{g_i},\cdots,\bar{g}_n)\in\widetilde{St}(V/L_i)
$$
for each $\langle g_i\rangle=L_i$ with $1\leq i\leq n$. Thus we get the formula
$$
\partial(0,g_1,\cdots,g_n)\mapsto \sum_{i=1}^n(-1)^i((L_i,V), \partial(0, \bar{g}_1, \cdots, \widehat{g_i},\cdots,\bar{g}_n)).
$$

\end{proof}
\

\subsubsection{The Formula of $d^1$ on Coefficients}

\begin{theorem}\label{d1coeff}
For $n\geq 2$, the formula of $d^1$ is given by
$$
\partial(0,g_1,\cdots,g_n)\longmapsto
$$
$$
-\sum_{i=1}^n(-1)^i((0, W_i), \partial(0, g_1,\cdots,\widehat{g_i},\cdots, g_n))
+\sum_{i=1}^n(-1)^i((L_i,V), \partial(0,\bar{g}_1, \cdots, \widehat{g_i},\cdots,\bar{g}_n)).
$$
\end{theorem}
\begin{proof}
It suffices to notice that, by our discussion in previous sections, there exists a commutative diagram
$$
\xymatrix{
H_{n-1}(\widetilde{\F_T}(d))\ar[rrr]^{d^1}\ar[d]_{\simeq} & & & \bigoplus_{c\to d\in\F_{T'}(d)}H_{n-2}(\widetilde{\F_U}(c))\ar[d]^{\simeq}\\
\widetilde{St}(d)\ar[rrr] & & & \left(\bigoplus_{(0,W)\in\F_{T'}(d)}\widetilde{St}(c)\right)\oplus\left(\bigoplus_{(L,V)\in\F_{T'}(d)}\widetilde{St}(V/L)\right)
}
$$
where $W=c\otimes K$.

\end{proof}
\

\begin{remark}\label{d1modular}
According to theorem \ref{compare5}, for $n\geq 3$, the formula for $d^1$ on coefficients theorem \ref{d1coeff} can be written as
$$
[g_1,\cdots,g_n]\longmapsto
$$
$$
-\sum_{i=1}^n(-1)^i((0, W_i), [g_1,\cdots,\widehat{g_i},\cdots, g_n])+\sum_{i=1}^n(-1)^i((L_i,V), [\bar{g}_1, \cdots, \widehat{g_i},\cdots,\bar{g}_n]).
$$
\end{remark}
\

\section{The Formula for $d^1|_{E^1_{n\geq 3,0}}$}

\subsection{The Induced Functor}

\begin{definition}\label{defind}
Suppose that $f :\C_1\to\C_2$ and $\F: \C_1\to\mathbf{Cat}$ are two functors and for each $y_2\in\C_2$ the category $\F_f(y_2)\neq\emptyset$. We define
$$
\Ind_f\F: \C_2\to \mathbf{Cat},\ \ \ \ \ \ y_2\mapsto \F_f(y_2)\int\F\circ\pi_{y_2}
$$
where $\pi_{y_2}:\F_f(y_2)\to\C_1$ is the projection.
\end{definition}
\begin{lemma}\label{lem2}
There is a canonical isomorphism
$$
H_*(\C_2, \Ind_f\F)\simeq H_*(\C_1, \F).
$$
\end{lemma}
\begin{proof}
We define the projection functor 
$$
p_1: \C_2\int \Ind_f\F\to \C_1\int\F,\ \ \ \ \ \{(y_1\to y_2, x)\ |\ x\in\F(y_1)\}\mapsto\{(y_1,x)\ |\ x\in\F(y_1)\}
$$
and
$$
s_1: \C_1\int\F\to  \C_2\int \Ind_f\F,\ \ \ \ \ \{(y_1,x)\ |\ x\in\F(y_1)\}\mapsto\{(y_1=y_1, x)\ |\ x\in\F(y_1)\}.
$$
It is easy to verify that $s_1$ is left adjoint to $p_1$ and hence the categories $\C_2\int \Ind_f\F$ and $\C_1\int\F$ are homotopy equivalent and hence the isomorphism we are looking for.

\end{proof}
\

This lemma tells us that we have a pair of canonical isomorphisms
$$
H_*(\C_1, \F)\simeq H_*(\C_2, \Ind_f\F).
$$
Combined with Eilenberg-Zilber-Cartier theorem, we obtain
$$
H_q(\C_1, H_*(\F))\simeq H_q(\C_2, H_*(\Ind(_f\F))),\ \ \ \ \ q\geq 0.
$$
\

Suppose that we are under the condition given at the beginning of section 5.1. We apply our discussion to 
$$
\xymatrix{
\C-\B\ar[rr]^{*}\ar[d]_{T'} & & \mathbf{Cat} \\
\D-\B\ar[rru]_{\Ind_{T'}*}
}
$$
Notice that $\Ind_{T'}*=\F_{T'}$, so the projection functor $(\D-\B)\int\F_{T'}\xrightarrow{p_1} C-\B$ given by $(c\to d)\to c$ admits a left adjoint $c\to(c=c)$ and hence is a homotopy equivalence. So we have 
$$
C_*\left((\D-\B)\int\F_{T'}, \widetilde{\F_{U_*}}\right)
$$
$$
=Ker\left(C_*\left((\D-\B)\int\F_{T'}, \F_{U_*}\right)\to C_*\left((\D-\B)\int\F_{T'}\right)\xrightarrow{\sim} C_*(\C-\B)\right)
$$
$$
=Ker\left(C_*\left((\D-\B)\int\F_{T'}, \F_{U_*}\right)\xrightarrow{p_{1*}} C_*(\C-\B, \F_U)\to C_*(\C-\B)\right).
$$
It follows that there is a canonical quasi-isomorphism 
$$
C_*\left((\D-\B)\int\F_{T'}, \widetilde{\F_{U_*}}\right)\xrightarrow{\sim}C_*(\C-\B, \widetilde{\F_U})
$$ 
and hence by Thomason's theorem (theorem \ref{thomason}) we get canonical isomorphisms 
$$
p_{1*}: H_q\left(\D-\B, H_{n-2}(\F_{T'}\int\widetilde{\F_{U_*}})\right)\simeq H_q\left((\D-\B)\int\F_{T'}, H_{n-2}(\widetilde{\F_{U_*}})\right)
$$
\begin{equation}\label{eq11}
\xrightarrow{\simeq}H_q(\C-\B, H_{n-2}(\widetilde{\F_U})).
\end{equation}
\

\subsection{The Formula for $d^1|_{E^1_{n\geq 3,0}}$}
Let us get back to the background set up in the previous section. By theorem \ref{compare5}, we will use modular symbols $[g_1,\cdots,g_n]$ instead of extended symbols in our formulae.  
\\

We notice that the canonical quasi-isomorphism $C_*(\D-\B, \F_{T'})\to C_*((\D-\B)\int\F_{T'})$ induced by Thomason's theorem \ref{thomason} sends the simple chain
$$
\xymatrix{
d_0\ar[r] & d_1\ar[r] & \cdots\ar[r] & d_p\\
c_0\ar[u]\ar@{=}[r] & c_0\ar@{=}[r] & \cdots\ar@{=}[r] & c_0
}
$$
to
$$
\xymatrix{
d_0\ar[r] & d_1\ar[r] & \cdots\ar[r] & d_p\\
c_0\ar[u]\ar@{=}[r] & c_0\ar[u]\ar@{=}[r] & \cdots\ar@{=}[r] & c_0\ar[u]
}
$$
\

Since $\F_{T'}(d)$ is discrete for $d\in\D-\C$ and if $d\in\C$ then $H_{n-1}(\widetilde{\F_T}(d))=0$. Apply the functor $C_q(\D-\B,-)$ to the formula obtained in proposition \ref{prop2}, we get
$$
C_q(\D-\B, H_{n-1}(\widetilde{\F_T}))\to C_q\left(\D-\B, H_{n-2}(\F_{T'}\int\widetilde{\F_{U_*}})\right)
$$
sending $0\neq\sum_d n_d(d_0\to\cdots\to d_{q}, \sum_gn_g[g_1,\cdots,g_n])$ with $d_j\in\D-\C$ for $0\leq j\leq q$ to
$$
\sum_d n_d\left(d_0\to\cdots\to d_{q}, \sum_gn_g\left(-\sum_{i=1}^n(-1)^i((0,W_i), [g_1,\cdots,\widehat{g_i},\cdots,g_n])\right.\right.
$$
$$
\left.\left.+((L_i,V), [\bar{g}_1,\cdots,\widehat{g_i},\cdots,\bar{g}_n])\right)\right).
$$
Combined with (\ref{eq11}), we have proved that
\begin{theorem} 
Let $A$ be an integral Noetherian ring and $Q^{\tf}(A)$ be Quillen's Q-construction over the category of finitely generated torsion-free modules. 
Then the differential
$E^1_{n,q}\xrightarrow{d^1_{n, q}} E^1_{n-1,q}$ is given by 
$$
\sum_d n_d\left(d_0\xrightarrow{f_1}\cdots\xrightarrow{f_q} d_{q}, \sum_gn_g[g_1,\cdots,g_n]\right)
$$
$$
\longmapsto \sum_{d,g} n_dn_g\left(\sum_i(-1)^i \left(-(c^0_i\to\cdots\to c^q_i, [g_1,\cdots,\widehat{g_i},\cdots, g_n])\right.\right.
$$
$$
\left.\left.+(c^{'0}_i\to\cdots\to c^{'q}_i, [\bar{g}_1,\cdots,\widehat{g_i},\cdots,\bar{g}_n])\right)\right).
$$
Here, $(0,W_i)=(0,c^0_i)$ is an admissible monomorphism and $c_i^j\to d_j=(f_j\circ\cdots\circ f_1)(0,W_i)$,  meanwhile $(L_i,V)=d_0\twoheadrightarrow c^{'0}_i$ is an admissible epimorphism and $(c^{'j}_i\to d_j)=(f_j\circ\cdots\circ f_1)(L_i,V)$.
\end{theorem}
\

\begin{remark}
For any $1\leq i\neq j\leq n$, the map $d^1$ sends the simple chain
$$
(d_0\to\cdots\to d_q, \sum_g n_g[g_1,\cdots,g_n])
$$
to
$$
\sum_i(-1)^i \left(-(c^0_i\to\cdots\to c^q_i, [g_1,\cdots,\widehat{g_i},\cdots, g_n])+(c^{'0}_i\to\cdots\to c^{'q}_i, [\bar{g}_1,\cdots,\widehat{g_i},\cdots,\bar{g}_n])\right).
$$
According to \cite[p17, diagram (4)]{Q1}, any morphism under Q-construction can be factored as an injective- followed-by-surjective map or a surjective-followed-by-injective map, and the factorizations are unique up to unique isomorphisms. In the language of proper layers, this can be written as
$$
(L_i,W_j)\circ(0,W_j)=(L_i,W_j)\circ(L_i,V)=(L_i,W_j).
$$
So if we apply $d^1$ again, the above chain will further be sent to
$$
\sum_{i=1}^{n}(-1)^i\left(\sum_{k=1}^{n-1}(-1)^k\left( (b^0_{ik}\to\cdots\to b^q_{ik}, [g_1,\cdots,\widehat{g_i},\cdots,\widehat{g_k},\cdots,g_{n}])\right.\right.
$$
$$
\left.\left.-(b^{'0}_{ik}\to\cdots\to b^{'q}_{ik}, [\bar{g}_1,\cdots,\widehat{g_i},\cdots,\widehat{g_k},\cdots,\bar{g}_{n}])\right)\right)+
$$
$$
\sum_{i=1}^{n}(-1)^i\left(\sum_{k=1}^{n-1}(-1)^k\left( -(b^{'0}_{ik}\to\cdots\to b^{'q}_{ik}, [\bar{g}_1,\cdots,\widehat{g_i},\cdots,\widehat{g_k},\cdots,\bar{g}_{n}])\right.\right.
$$
$$
\left.\left.+(b^{''0}_{ik}\to\cdots\to b^{''q}_{ik}, [\bar{g}_1,\cdots,\widehat{g_i},\cdots,\widehat{g_k},\cdots,\bar{g}_{n}])\right)\right)
$$
such that $k=j$ if $1\leq j<i$ and $k=j-1$ if $i<j\leq n$. Here, $b^0_{ik}=\langle g_1,\cdots,\widehat{g_i},\cdots\widehat{g_k},\cdots,g_{n}\rangle\cap c^0_i$, $b^{'0}_{ik}\to c^0_i=(L_k,W_i)$, $b^{'0}_{ik}\to c^{'0}_i=(L_i,W_k)$ and $b^{''0}_{ik}\to c'_0=(L_k, V/L_i)=(L_k\oplus L_i, V)$. 
Thus we find that in the above formula, 
a summand with same two entries omitted appears twice with opposite signs, so the resulting alternating sum is zero which verifies that $d^1\circ d^1=0$. 
\end{remark}




\cleardoublepage 
\renewcommand{\indexname}{Index des notations} 
\printindex 

\end{document}